\documentclass[11pt]{amsart}

\usepackage{amsmath,amsthm,indentfirst}
\usepackage{amssymb}
\usepackage{amsfonts}
\usepackage{amscd}
\usepackage[latin1]{inputenc}
\usepackage{hyperref}
\usepackage{ifthen, amsfonts, amssymb, graphicx, srcltx, mathrsfs, xfrac,
amsmath}
\usepackage{enumerate}

\theoremstyle{plain}
\newtheorem*{maintheorem*}{Main Theorem}

\newtheorem*{thmnarrowend}{Theorem 1.1}
\newtheorem*{thmdivsymb}{Theorem 6.7}

\newtheorem*{remark*}{Remark}
\newtheorem*{conjecture*}{Conjecture}
\newtheorem*{prop*}{Proposition}
\newtheorem{thm}{Theorem}[section]
\newtheorem{cor}[thm]{Corollary}
\newtheorem{lem}[thm]{Lemma}
\newtheorem{prop}[thm]{Proposition}

\theoremstyle{definition}

\newtheorem*{proofc*}{Proof of Theorem C}

\newtheorem{definition}[thm]{Definition}

\newtheorem{remark}[thm]{Remark}
\newtheorem{notation}[thm]{Notation}
\newtheorem{claim}[thm]{Claim}

\DeclareMathOperator{\CAT}{CAT}
\DeclareMathOperator{\Teich}{Teich}
\DeclareMathOperator{\Mod}{Mod}

\DeclareMathOperator{\grad}{grad}
\DeclareMathOperator{\WP}{WP}
\DeclareMathOperator{\base}{base}
\DeclareMathOperator{\tw}{tw}
\DeclareMathOperator{\length}{length}
\DeclareMathOperator{\inj}{inj}

\DeclareMathOperator{\diam}{diam}

\DeclareMathOperator{\spn}{span}

\numberwithin{equation}{section}

\begin{document}


\title[Asymptotics and divergence of Weil-Petersson geodesics]{Asymptotics of a class of Weil-Petersson geodesics and divergence of Weil-Petersson geodesics}

\date{\today}


\author{Babak Modami}
\address{Department of Mathematics, University of Illinois at Urbana-Champaign, 1409 W Green ST, Urbana, IL 61801}
\email{bmodami@illinois.edu}
\thanks{The author was partially supported by NSF grant DMS-1005973.}


\date{\today}

\maketitle
\begin{abstract}
We show that the strong asymptotic class of Weil-Petersson (WP) geodesic rays with narrow end invariant and bounded annular coefficients is determined by the forward ending laminations of the geodesic rays. This generalizes the Recurrent Ending Lamination Theorem of Brock, Masur and Minsky. As an application we provide a symbolic condition for divergence of WP geodesic rays in the moduli space.
\end{abstract}

\setcounter{tocdepth}{3}
\tableofcontents


\section{Introduction}

The Weil-Petersson (WP) metric is a Riemannian metric on the moduli space of Riemann surfaces. Over the last decade various aspects of the geometry and dynamics of the metric have been studied, however in comparison with the Teichm\"{u}ller metric-- the most well studied metric on the moduli space-- much less is known about this metric. The WP metric is incomplete with sectional curvatures asymptotic to $0$ and $-\infty$ in the completion.  These features, in particular, prevent applying some of the standard techniques to study the global geometry and dynamics of the metric. For example, the Shadow Lemma for construction geodesic rays with specific behavior. WP geodesic rays are not necessarily visible \cite{brockcommunication}. The Weil-Petersson geodesic flow is not uniformly hyperbolic. Also, curently Markov partitions are not available for coding of the WP geodesic flow.
\medskip

Brock, Masur and Minsky \cite{bmm1}, in analogy with the vertical geodesic lamination of a Teichm\"{u}ller geodesic, introduced a notion of ending lamination for Weil-Petersson geodesic rays. They conjectured that ending laminations, or a modification of them, can be used to parametrize the visual boundary of the WP metric and also the stable and unstable foliations of the WP geodesic flow. Furthermore, it is conjectured that end invariants and the associated subsurface coefficients provide a kind of symbolic coding for Weil-Petersson geodesics in the moduli space. 
\medskip

Further, Brock, Masur and Minsky \cite{bmm1}, \cite{bmm2} explored several aspects of the mentioned conjectures. Significantly, they proved that the forward ending lamination determines the strong asymptotic class of recurrent WP geodesic rays to the thick part of the moduli space. Moreover, they showed that the bounded combinatorics of end invariants is equivalent to co-boundedness of the geodesic; the geodesic projects to a compact subset of the moduli space. These results have dynamical consequences, among which are the topological transitivity of the Weil-Petersson geodesic flow on the moduli space and unboundedness the topological entropy of the WP flow.
\medskip

In \cite{wpbehavior} we considered WP geodesics with narrow end invariants, end invariants with a certain constraint on subsurfaces with a big subsurface coefficient (see Definition \ref{def : narrow}), and constructed examples of closed WP geodesics in the thin parts of the moduli space as well as divergent WP geodesic rays with minimal filling ending laminations.
 \medskip

In this paper we show that the strong asymptotic class of a WP geodesic ray with narrow end invariant and bounded annular coefficients is determined by the forward ending lamination.

\begin{thm}\label{thm : strongasym}\textnormal{(Narrow ending lamination theorem)}
The strong asymptotic class of a WP geodesic ray with narrow end invariant and bounded annular coefficients is determined by the forward ending lamination. 
\end{thm}
The strong asymptotic class of a geodesic ray $r$ is the set of all the rays $r'$ with $d(r(t),r'(t))\to 0$ as $t\to \infty$. The class of WP geodesic rays with narrow ending invariant and bounded annular coefficients contains geodesic rays which are not recurrent to any compact subset of the moduli space (divergent rays); see \cite[\S8]{wpbehavior}. Heuristically these geodesic rays avoid all asymptotic flats in the WP metric and exhibit features of geodesics in manifolds with negative sectional curvatures which are bounded away from $0$. This theorem is a generalization of the following result from \cite{bmm1}.

\begin{thm}\textnormal{(Recurrent ending lamination theorem)}\label{thm : recel}
The strong asymptotic class of a WP geodesic ray recurrent to a compact subset of the moduli space is determined by its forward ending lamination. 
\end{thm}

These theorems address the parametrization of the visual boundary of the WP metric and characterization of the stable and unstable foliations of the WP geodesic flow using laminations.
\medskip

For the proof of Theorem \ref{thm : strongasym} we use the control of the length-functions along WP geodesics developed in \cite{wpbehavior} and ruled surfaces as in \cite{bmm1}. The new ingredient here is the strict uniform contraction property of the nearest point projection to WP geodesic segments close to the thick part of a stratum which is not the product of lower complexity strata; see $\S$ 5, in particular, Theorems \ref{thm : contnpt} and \ref{thm : normJbd}. The contraction property is proved using some of Wolpert's estimates on the WP metric, the WP Levi-Civita covariant derivatives and sectional curvatures in the thin part of the Teichm\"{u}ller space and compactness arguments. 

For WP geodesic rays with prescribed itinerary, as in Theorem \ref{thm : inftyray}, using the contraction property we can guarantee the existence of regions with a definite negative total curvature on ruled surfaces with one side on the ray with prescribed itinerary; see $\S 4$ and the proof of Theorem \ref{thm : asymsasym}.

A geodesic ray $r$ in a metric space is visible if for any other geodesic ray $r'$ there is an infinite geodesic (strongly) asymptotic to $r$ in the forward time and (strongly) asymptotic to $r'$ in the backward time. In a complete Riemannian manifold with negative sectional curvatures bounded away from $0$ every geodesic ray is visible. For the notion of visibility and some of its dynamical consequences see \cite{ebergoedflownegcurv}. In the regions with a definite negative total curvature as above we are able to pick up enough negative curvature on the ruled surface so that using a variation of the Gauss-Bonnet formula asymptotic convergence to the ray with prescribed itinerary is guaranteed (Theorem \ref{thm : asymsasym}). Using a similar technique, we overcome the difficulty caused by the fact that the sectional curvatures of the WP metric are not bounded away from $0$ in the thin part of the Teichm\"{u}ller space, and we prove visibility of the class of geodesic rays with narrow end invariant and bounded annular coefficients in Theorem \ref{thm : narrowvis}.

Finally, as an application of our ending lamination theorem we prove a symbolic condition in terms of subsurface coefficients for divergence of WP geodesic rays in the moduli space:
\begin{thmdivsymb} \textnormal{(Divergence condition)}
Given $A,R,R'>0$. Let $(\nu^{-},\nu^{+})$ be an $A-$narrow pair on a Riemann surface $S$ with $R'-$bounded annular coefficients, and suppose that $d_{S}(\nu^{-},\nu^{+})\leq R$. Then a WP geodesic ray with end invariant $(\nu^{-},\nu^{+})$ is divergent in the moduli space $\mathcal{M}(S)$.
\end{thmdivsymb}
\medskip

\textbf{Acknowledgement:} I am so grateful to Yair Minsky for many invaluable discussions in the course of this work. I would also like to thank Scott Wolpert for several communications related to this work.

\section{Background}

\begin{notation} Let $f,g:X\to \mathbb{R}^{\geq 0}$ be two functions. Let $K\geq 1$ and $C\geq 0$ be two
 constants. We write $f\asymp_{K,C} g$ if
$$\frac{1}{K}g(x)-C\leq f(x)\leq Kg(x)+C,$$
 holds for every $x\in X$.
\end{notation}

\subsection{Curve complexes and hierarchy paths} \label{subsec : ccplx}Let $S=S_{g,b}$ be a finite type, closed, orientable surface with genus $g$ and $b$ punctures or boundary components. Define the complexity of $S$ by $\xi(S)=3g-3+b$. A subsurface of $S$ is an embedded, closed subsurface of $S$ with non-perpheral boundary curves.

The \textit{curve complex} of the surface $S$, denoted by $\mathcal{C}(S)$, is a flag complex which serves to organize isotopy classes of simple closed curves on a surface. The complex is defined as follows: When $\xi(S)>1$, each vertex of the complex is the isotopy class of an essential simple closed curve, with an edge between each pair of isotopy classes with disjoint representatives on $S$. In the same fashion there is a $k$ simplex corresponding to any set of $k+1$ pair-wise simple closed curves on $S$. When $\xi(S)=1$, $S$ is $S_{0,4}$ or $S_{1,1}$. The definition of the curve complex is the same, except that there is an edge between any pair of isotopy classes of curves with representatives with intersection number $1$ in the case of $S_{1,1}$ and $2$ in the case of $S_{0,4}$.  

When $Y\subset S$ is an annulus with essential core curve the definition is slightly different. Let $\widetilde{Y}$ be the annular cover of $S$ to which $Y$ lifts homeomorphically. There is a natural compactification of $\widetilde{Y}$ to a closed annulus $\widehat{Y}$ which is obtained from the compactification of the Poincare disk (the universal cover of $S$) by the closed disk. A vertex of $\mathcal{C}(Y )$ is associated to an arc connecting the two boundary components of $\widehat{Y}$ modulo isotopies that fix the endpoints (isotopy classes of arcs relative to the boundary). There is an edge between two vertices which have representatives with disjoint interiors.  

We equip the curve complex with a distance by making each simplex Euclidean with side lengths 1, and denote the distance by $d_{S}=d_{\mathcal{C}(S)}$. One can easily verify that the curve complex of any annular subsurface is quasi-isometric to $\mathbb{Z}$. Morever Masur and Minsky in their seminal work \cite{mm1} proved that the curve complex of $S$ is $\delta-$hyperbolic with $\delta$ depending only on the topological type of $S$. 

A {\it multi-curve} on the surface $S$ is a collection of pair wise disjoint curves. Let $\sigma$ and $\tau$ be two multi-curves. We say that $\sigma$ and $\tau$ {\it overlap} and write $\sigma\pitchfork\tau$, if there are curves $\alpha\in\sigma$ and $\beta\in\tau$ that intersect each other essentially (i.e. can not be realized as disjoint curves on $S$). 
\medskip

\noindent{\bf Laminations:} Fix a complete finite area hyperbolic metric on $S$. A\textit{ geodesic lamination} on $S$ is a closed subset of $S$ consisting of complete, simple geodesics. In particular, $\mathcal{C}_{0}(S)\subset \mathcal{GL}(S)$. Geodesic laminations provide a natural completion for the curve complex and Teichm\"{u}ller space. 

Each geodesic lamination can be equipped with a {\it transverse measure} \cite{phtraintr}. The pair of a lamination and a transverse measure is called a {\it measured geodesic lamination}. We denote the space of measured geodesic laminations equipped with the weak$^{*}$ topology by $\mathcal{ML}(S)$. $\mathbb{R}^{+}$ acts on $\mathcal{ML}(S)$ by rescaling measures, each equivalence class is called a {\it projective measured lamination}. We denote the quotient space by $\mathcal{PML}(S)$.

Recall that the curve complex of $S$ is Gromov hyperbolic. By a result of Klarreich \cite{bdrycc} the Gromov boundary of the curve complex is identified with the {\it ending laminations space} $\mathcal{EL}(S)$. The space $\mathcal{EL}(S)$ is the image of projective measured laminations with minimal filling support in $\mathcal{GL}(S)$ under the measure forgetful map. Moreover, $\mathcal{EL}(S)$ is equipped with the topology induced from the topology of $\mathcal{PML}(S)$ via the measure forgetful map.
\medskip

 \noindent{\bf Pants decomposition and markings:} A {\it pants decomposition} $P$ on $S$ is a maximal set of pair wise disjoint curves on $S$. A {\it (partial) marking} $\mu$ is obtained from a pants decomposition $P$ by adding transversal curves to (some) all of the curve in the pants decomposition. We call $P$ the base of $\mu$ and denote it by $\base(\mu)$. The set of all pants decompositions can be turned into a metric graph which is called the {\it pants graph}. For this purpose we put a length one edge between any two pants decompositions which differ by an elementary move. Similarly the markings can be turned into a metric graph called the {\it marking  graph}. For more detail see \cite{mm2}. 
\medskip

\noindent\textbf{Subsurface coefficient:} 
An {\it essential subsurface} $Y\subseteq S$ is a compact, connected subsurface of $S$ whose boundary consists of essential curves in $S$ or boundary curves of $S$ which is not a $3-$holed sphere. In this paper we do not distinguish between a subsurface and its isotopy class.

 Let $Y\subseteq S$ be an essential subsurface. We define the \textit{subsurface projection map} 
 $$\pi_{Y}:\mathcal{GL}(S)\to \mathcal{PC}_{0}(Y)$$
where $\mathcal{PC}_{0}(Y)$ is the power set of $\mathcal{C}_{0}(Y)$ as follows: Equip $S$ with a complete hyperbolic metric and realize all curves and laminations geodesically. Let $\lambda\in \mathcal{GL}(S)$. Suppose that $Y$ is a non-annular subsurface. If $\lambda\cap Y=\emptyset$, then define $\pi_{Y}(\lambda)=\emptyset$. Otherwise, consider the set of arcs with end points on $\partial{Y}$ or at cusps in the intersection locus $\lambda\cap Y$. Identify all of the curves and arcs that are isotopic. Through these isotopies the end points of arcs are allowed to move within $\partial{Y}$. The projection $\pi_{Y}(\lambda)$ is the union of the boundary curves of a regular neighborhood of $a\cup\partial{Y}$ where $a$ is an arc we obtained above, and all the closed curves we obtained above. 

For an annulus $Y$ the subsurface projection is the set of component arcs of the lift of $\lambda$ to $\widehat{Y}$ (the compactification of the annular cover of $S$ to which $Y$ lifts homeomorphically) that connect the two boundaries of $\widehat{Y}$. For more detail see \cite[\S 2]{mm2}.
 
 The projection of a multi-curve $\sigma$ to a subsurface $Y$ is the union of the projections $\pi_{Y}(\alpha)$ of all $\alpha\in\sigma$. Let $\mu$ be a partial marking. If $Y$ is an annular subsurface with core curve $\alpha\in\base(\mu)$, then $\pi_{Y}(\mu)$ is the set of transversal curves of $\alpha$. Otherwise, $\pi_{Y}(\mu)=\pi_{Y}(\base(\mu))$.
 \medskip
 
  The {\it $Y$ subsurface coefficient} of two laminations or (partial) markings $\mu$ and $\mu'$ is defined by 
  \begin{equation}\label{eq : subsurfcoeff}d_{Y}(\mu,\mu')=\min\{d_{Y}(\gamma,\gamma') : \gamma\in \pi_{Y}(\mu), \gamma'\in \pi_{Y}(\mu')\}.\end{equation}
   Note that the subsurface coefficients are an analogue of continued fraction expansions of real numbers. We use subsurface coefficients to study geodesics on moduli spaces, similar to the role of continued fraction expansion in the study of geodesics on the modular surface.
\medskip

We denote by $\diam_{Y}(\mu)$ the diameter of $\pi_{Y}(\mu)$ viewed as a subset of $\mathcal{C}(Y)$. The following lemma is a straightforward consequence of \cite[Lemma 2.3]{mm2}.

 \begin{lem}\label{lem : diamproj} 
 Let $\mu$ be a (partial) marking on a surface $S$. For any essential subsurface $Y\subseteq S$ we have
 $$\diam_{Y}(\mu)\leq 2.$$
\end{lem}
\medskip

The {\it hierarchy paths} introduced by Masur and Minsky \cite{mm2} comprise a transitive family of quasi-geodesics in the pants graph of a surface with quantifiers depending only the topological type of the surface. These quasi-geodesics are constructed from hierarchies of geodesics in subsurfaces of the surface. The main feature of hierarchy paths is that their properties are encoded in their end points and the associated subsurface coefficients. For a list of the properties of hierarchy paths see \cite{bmm2}, \cite{wpbehavior}. 

Let $(\mu^{-},\mu^{+})$ be a pair of (partial) markings or laminations. Let $\rho:[m,n]\to P(S)$, where $[m,n]\subseteq \mathbb{Z}\cup \{\pm \infty\}$, be a hierarchy path with $\rho(m)=\mu^{-}$ and $\rho(n)=\mu^{+}$. An important property of the hierarchy path is the following: There are subsurfaces $Y$ called component domains; corresponding to each component domain there is a connected subinterval $J_{Y}\subset [m,n]$ such that $\partial{Y}\subset \rho(i)$ for all $i\in J_{Y}$. Moreover there is a constant $M_{1}>0$ depending only on the topological type of the surface, so that any subsurface $Y$ with $d_{Y}(\mu^{-},\mu^{+})>M_{1}$ is a component domain.
 \medskip

Using the machinery of hierarchies  Masur and Minsky \cite{mm2} established the following quasi-distance formula: 

Given $A\geq M_{1}$, there exist constants $K\geq 1$ and $C\geq 0$ such that for any $P,Q\in P(S)$ we have
\begin{equation}\label{eq : dsf}d(P,Q)\asymp_{K,C}\sum_{\substack{Y\subseteq S\\ \text{non-annular}}}\{d_{Y}(P,Q)\}_{A} \end{equation}
Here the cut off function $\{.\}_{A}:\mathbb{R}\to \mathbb{R}^{+}$ is defined by $\{a\}_{A} = \Big\lbrace 
  \begin{array}{cl}
  a & \text{if}\; a \geq A\\
  0 & \text{if}\; a < A \\
  \end{array}.$
We call $A$ the threshold constant and $K$ and $C$ the corresponding multiplicative and additive constants.
\medskip

In the following we define two constraints on subsurface coefficients which would be used in this paper. 
\begin{definition}(\textnormal{Narrow pair})\label{def : narrow}
Let $A>0$. An $A-$narrow pair of (partial) markings or laminations $(\mu^{-},\mu^{+})$ is a pair such that for an essential subsurface $Z$ if 
$$d_{Z}(\mu^{-},\mu^{+})>A,$$
 then $Z$ is large, i.e. each connected component of $S\backslash Z$ is an annulus or a three holed sphere. In \cite{wpbehavior} we proved that any hierarchy path between a narrow pair is stable in the pants graph. 
\end{definition}

\begin{definition}\textnormal{(Bounded combinatorics})\label{def : bddcomb}
 Given $R,R'>0$. Let $\rho:[m,n]\to P(S)$ be a hierarchy path. Let $i,j\in [m,n]$. We say that $\rho$ has non-annular $(R,R')-$bounded combinatorics over $[i,j]\subset J_{Z}$ if
$$d_{Y}(\rho(i),\rho(j))\leq R$$
for every essential non-annular subsurface $Y\subseteq Z$. Moreover,
$$d_{\gamma}(\rho(i),\rho(j))\leq R',$$
for every $\gamma\in\mathcal{C}_{0}(Z)$.
\end{definition}

\subsection{Weil-Petersson metric}\label{subsec : wp}
Let $S$ be a surface with only punctures (no boundary curves). A point in the {\it Teichm\"{u}ller space of the surface $S$}, denote by $\Teich(S)$, is a complete finite area hyperbolic surface equipped with a diffeomorphism $h:S\to x$. The diffeomorphism $h$ is a marking of $x$. Two marked surfaces $h_{1}:S\to x_{1}$ and $h_{2}:S\to x_{2}$ define the same point in $\Teich(S)$ if and only if $h_{2} \circ h_{1}^{-1}:x_{1}\to x_{2}$ is isotopic to an isometry. The {\it moduli space of $S$}, denoted by $\mathcal{M}(S)$, is the quotient of $\Teich(S)$ by the action of the mapping class group of $S$, denoted by $\Mod(S)$, on $\Teich(S)$ by precomposition (remarking). 

The $\epsilon-$thick part of the Teichm\"{u}ller space consists of all $x\in\Teich(S)$ with injectivity radius $\inj(x)> \epsilon$. The $\epsilon-$thin part consists of all $x$ with $\inj(x)<\epsilon$. Suppose that $\epsilon$ is small enough that by the Collar Lemma (see \cite[\S 4.1]{buser}) there is no pair of intersecting closed geodesics of length less than or equal to $\epsilon$ on a complete hyperbolic surface. Given a multi-curve $\sigma$ we define the open regions in Teichm\"{u}ller space
\begin{itemize}
\item$ U_{\epsilon}(\sigma):=\{x\in \Teich(S) : \ell_{\alpha}(x)<\epsilon \;\text{for every}\; \alpha\in \sigma\}$,
\item $U_{\epsilon,\epsilon'}(\sigma):=\{x\in \Teich(S) : \ell_{\alpha}(x)<\epsilon \;\text{for every}\; \alpha\in \sigma \;\text{and}\; \ell_{\gamma}(x)> \epsilon' \;\text{for every}\; \gamma\notin \sigma\}$.
\end{itemize}
For a comprehensive introduction to the Weil-Petersson metric and its properties we refer the reader to \cite{wol},\cite{wolb}. Here we only recall some of the properties of the metric which are important for us. The Weil-Petersson (WP) metric is a Riemannian metric with negative sectional curvatures on the Teichm\"{u}ller space. The metric is incomplete due to possibility of pinching curves along paths with finite WP length. However there is a unique WP geodesic between any two points in the Teichm\"{u}ller space so the metric is geodesically convex. The WP sectional curvatures are negative and asymptotic to both $0$ and $-\infty$ asymptotic to completion. The completion is a $\CAT(0)$ metric. The metric is invariant under the action of $\Mod(S)$ and descends under the natural orbifold cover to a metric on $\mathcal{M}(S)$. 

 The Weil-Petersson completion of the Teichm\"{u}ller space $\overline{\Teich(S)}$ is the disjoint union of $\sigma-$strata denoted by $\mathcal{S}(\sigma)$ where $\sigma$ is a multi-curve. $\mathcal{S}(\sigma)$ consists of nodal Riemann surfaces at $\sigma$. Equivalently each point in $\mathcal{S}(\sigma)$ is a marked complete hyperbolic surface with a pair of cusps for each curve in $\sigma$. The WP completion of the Teichm\"{u}ller space descends to the Deligne-Mumford compactification of the moduli space $\overline{\mathcal{M}(S)}$ via the action of the mapping class group. $\overline{\mathcal{M}(S)}$ is the union of finitely many strata. Each stratum is the quotient of the strata of $\overline{\Teich(S)}$ which are identified by the natural extension of the action of $\Mod(S)$ to the completion.

The WP metric has the following {\it non-refraction property}: Given points $x,y\in\overline{\Teich(S)}$, the interior of the unique geodesic connecting $x$ and $y$ is in the smallest stratum (with respect to inclusion of strata) that contains $x$ and $y$. The following strengthened version of Wolpert's Geodesic Limit Theorem (see \cite{wols} and \cite{bmm2}) proved in $\S 4$ of \cite{wpbehavior} provides a limiting picture for a sequence of bounded length WP 
geodesic segments in the Teichm\"{u}ller space. We need this result for compactness arguments in $\S 5$.
\begin{thm}  \textnormal{(Geodesic limits)} \label{thm : geodlimit} 
Given $T>0$, let $\zeta_{n}:[0,T]\to \overline{\Teich(S)}$ be a sequence of WP
geodesic segments parametrized by arc-length.  
After possibly passing to a subsequence, there are a partition
$0=t_{0}<....<t_{k+1}=T$ of $[0,T]$, possibly empty multi-curves
$\sigma_{0},...,\sigma_{k+1}$, a multi-curve $\hat{\tau}\equiv
\sigma_{i}\cap \sigma_{i+1}$ for $i=0,1,...,k$, and a piece-wise
geodesic 
 $$\hat{\zeta}:[0,T]\to \overline{\Teich(S)},$$
 with the following properties
\begin{enumerate}[(1)]
\item\label{gl : tau} $\hat{\zeta}((t_{i},t_{i+1}))\subset \mathcal{S}(\hat{\tau})$ for $i=0,...,k$, 
\item \label{gl : sigma}$\hat{\zeta}(t_{i})\in \mathcal{S}(\sigma_{i})$ for $i=1,...,k$,

Given a multi-curve $\sigma$ denote by $\tw(\sigma)$ the subgroup of $\Mod(S)$ generated by positive Dehn twists about the curves in $\sigma$. There are elements of the mapping class group $\psi_{n}$ for each $n\in\mathbb{N}$ so that either $\psi_{n}\equiv I$ or $\psi_{n}$ is an unbounded sequence, and $\mathcal{T}_{i,n}\in \tw(\sigma_{i}-\hat{\tau})$ for $ i=1,...,k+1$ and any $n\in\mathbb{N}$ so that
\item\label{gl : converge} For any $t\in [0,t_{1}]$, $\psi_{n}(\zeta_{n}(t))\to \hat{\zeta}(t)$ as $n\to \infty$. For each $i=1,...,k$ and each $n\in\mathbb{N}$, let $\varphi_{i,n}=\mathcal{T}_{i,n}\circ ...\circ \mathcal{T}_{1,n}\circ \psi_{n}$. Then for any $t\in[t_{i},t_{i+1}]$, $\varphi_{i,n}(\zeta_{n}(t))\to \hat{\zeta}(t)$, as $n\to \infty$.
 \end{enumerate}
\end{thm}

\subsection{End invariants}\label{subsec : endinv}
There exists a constant $L_{S}$ (Bers constant) depending only on the topological type of the surface $S$ with the property that any complete finite area hyperbolic metric on $S$ has a pants decomposition (Bers pants decomposition) such that the length of every curve in $P$ is at most $L_{S}$; see \cite[\S 5]{buser}. A Bers curve is a curve in a Bers pants decomposition. A {\it Bers marking} is a (partial) marking obtained from a Bers pants decomposition by adding transversal curves with representatives of minimal length. Given a point $x\in\overline{\Teich(S)}$, suppose that $x\in\mathcal{S}(\sigma)$. A Bers pants decomposition of $x$, denoted by $Q(x)$, is the union of Bers pants decompositions of the connected components of $S\backslash\sigma$ and the multi-curve $\sigma$. A Bers marking of $x$, denoted by $\mu(x)$, is obtained from $Q(x)$ with no transversal for 
curves in $\sigma$. 
\medskip

By Brock's Quasi-Isometry Theorem \cite{br} the coarse map 
$$Q:\overline{\Teich(S)}\to P(S)$$
 which assigns to $x\in \Teich(S)$ a Bers pants decomposition of $x$ is a quasi-isometry with constants $K_{\WP}\geq 1$ and $C_{\WP}\geq 0$ depending only on the topological type of $S$. 
 
\begin{definition}\textnormal{(Ending measured lamination)}\label{def : endmeas}
The  weak$^{*}$ limit in $\mathcal{ML}(S)$ of any infinite sequence of weighted, distinct Bers curves along a WP geodesic ray $r$ is an ending measured lamination of $r$. 
\end{definition}

For any $\alpha\in\mathcal{C}_{0}(S)$, the {\it $\alpha-$length-function} assigns to a point $x\in\Teich(S)$ the length of the geodesic representative of $\alpha$ on the hyperbolic surface $x$. This notion of length-function has a natural extension to the space of measured laminations $\mathcal{ML}(S)$. Let $\mathcal{L}\in\mathcal{ML}(S)$. We denote the value of the $\mathcal{L}-$length-function at a point $x$ by $\mathcal{L}$ by $\ell_{\mathcal{L}}(x)$.

The convexity of length-functions along WP geodesics proved by Wolpert (see e.g. \cite[\S 3]{wol}) asserts that
\begin{thm}\label{thm : convlf}
Let $g$ be a WP geodesic. For any $\alpha\in\mathcal{C}_{0}(S)$, $\ell_{\alpha}(g(t))$ is a convex function. Similarly, for any $\mathcal{L}\in\mathcal{ML}(S)$, $\ell_{\mathcal{L}}(g(t))$ is a convex function. 
\end{thm}

In \cite{bmm2} the following notion of \textit{ending lamination} for WP geodesic rays is introduced. Its existence relies on the convexity of length-functions along WP geodesics and properties of $\CAT(0)$ spaces. Let $r:[0,a)\to\overline{\Teich(S)}$ be a WP geodesic ray.

\begin{definition}\textnormal{(Ending Lamination)}\label{def : endlam}
The union of pinching curves along a WP geodesic ray $r$ and the geodesic laminations arising as supports of all ending measured laminations of $r$ is the ending lamination of $r$. Where a pinching curve of $r$ is a curve $\alpha$ such that $\ell_{\alpha}(r(t))\to 0$ as $t\to a$.
\end{definition}

\begin{definition}\textnormal{(End invariant of Weil-Petersson geodesics)} To each open end of a geodesic $g:(a,b)\to\overline{\Teich(S)}$ (we assume that $0\in(a,b)$) we associate an end invariant which is a partial marking or a lamination. If the forward trajectory $g|_{[0,b)}$ can be extended to $b$ such that $g(b) \in \overline{ \Teich(S)}$ then the forward end invariant $\nu^{+}(g)$ is any Bers marking $\mu(g(b))$ (there are finitely many of them). Otherwise, $\nu^{+}(g)$ is the ending lamination of the forward trajectory ray $g|_{[0,b)}$ which was defined above. We define the backward end invariant $\nu^{-}(g)$ similarly by considering the backward trajectory $g|_{[0,a)}$. We call the pair $(\nu^{-},\nu^{+})$ the end invariant of $g$.
\end{definition}
We recall two important properties of  the ending measured laminations proved in \cite[\S 2]{bmm1}.
\begin{lem}\textnormal{(Decreasing of length along WP geodesic rays)}\label{lem : dlendmeas}
Let $\mathcal{L}$ be any ending measured lamination of a WP geodesic ray $r$, then $\ell_{\mathcal{L}}(r(t))$ is a decreasing function.
\end{lem}
\begin{lem} \label{lem : limraylimmeas} Let $r_{n} \to r$ be a convergent sequence of rays in the WP visual sphere at $x$. Then if $\mathcal{L}_{n}$ is any sequence of ending measured laminations or weighted pinching curves for $r_{n}$, any representative $\mathcal{L} \in \mathcal{ML} (S)$ of the limit of the projective classes $[\mathcal{L}_{n}]$ in $\mathcal{PML}(S)$ has bounded length along the ray $r$.
\end{lem}

Let $\nu^{+}$ be a measurable geodesic lamination. Suppose that there is a collection of pairwise disjoint subsurfaces $Z_{a}$, $a=1,...,m$, with $\xi(Z_{a})\geq 1$ so that any simple closed curve in $S\backslash \bigcup_{a=1}^{m}Z_{a}$ is isotopic to a boundary curve of one of the subsurfaces $Z_{a}$, and moreover that $\nu_{a}$ the restriction of $\nu^{+}$ to $Z_{a}$ is minimal and fills $Z_{a}$. For each $a=1,...,m$ let $\mathcal{L}_{a}$ be a measured lamination supported on $\nu_{a}$. Let $\gamma^{a}_{n}\in \mathcal{C}(Z_{a})$ be a sequence of curves so that the projective classes $[\gamma^{a}_{n}]$ converge to $[\mathcal{L}_{a}]$ in $\mathcal{PML}(Z_{a})$  as $n\to \infty$. For each $n\in\mathbb{N}$, let $Q_{n}$ be a pants decomposition that contains $\{\partial{Z}_{a},\gamma^{a}_{n}\}_{a=1}^{m}$. Let $c_{n}$ be the maximally nodal hyperbolic surface at $Q_{n}$. Let $[x,c_{n}]$ be the WP geodesic segment connecting a base point $x$ in the interior of the Teichm\"{u}ller space to $c_{n}$. Denote the parametrization of $[x,c_{n}]$ by arc-length by $r_{n}$. In \cite[\S 8]{wpbehavior} we proved

\begin{lem}\textnormal{(Infinite ray)}\label{lem : inftyray}
After possibly passing to a subsequence the geodesic segments $r_{n}$ converge to an infinite ray $r$ in the visual sphere of the WP metric at $x$. Moreover, the length of each measured lamination $\mathcal{L}_{a}$, $a=1,...,m$, and each curve $\alpha\in\partial{Z}_{a}$ is bounded along $r$.
\end{lem}

Let $(\nu^{-},\nu^{+})$ be a narrow pair. The narrow condition implies that there is at most one subsurface $Y$ with $\xi(Y)\geq 1$ so that the restriction of $\nu^{+}$ to $Y$ is minimal and fills $Y$. Suppose that such a component exists. Let $\rho$ be a hierarchy between $\nu^{-}$ and $\nu^{+}$. There is an $N\in\mathbb{N}$ so that for all $n\in\mathbb{N}$ we have
$$\partial{Y}\subset\rho(n+N).$$
Note that if $Y=S$, then this statement vacuously holds. For each $n\in\mathbb{N}$, let $Q_{n}=\rho(n+N)$ and $\gamma_{n}$ be a curve in $\pi_{Y}(Q_{n})$. After possibly passing to a subsequence the projective classes of $\gamma_{n}$ converge to the projective class of a measured geodesic lamination supported on $\nu'$. Let $x\in\Teich(S)$ be a point with a Bers marking $\nu^{-}$. Let $c_{n}$ be the maximally nodal hyperbolic surface at $Q_{n}$. As before, let $r$ be the limit of the geodesic segments $[x,c_{n}]$ after possibly passing to a subsequence. In \cite[\S 8]{wpbehavior} we proved:

\begin{lem}\label{lem : inftyraynu}
The forward ending lamination of $r$ contains $\nu'$.
\end{lem}

\section{Combinatorial control}

Let $g:[a,b]\to \Teich(S)$ be a WP geodesic with $A-$narrow end invariant $(\nu^{-},\nu^{+})$. Let $\rho:[m,n]\to P(S)$ be a hierarchy path between $\nu^{-}$ and $\nu^{+}$. In \cite[\S 5]{wpbehavior} we proved that a hierarchy path with narrow end points is $d-$stable. Where $d:\mathbb{R}^{\geq 1}\times \mathbb{R}^{\geq 0}\to \mathbb{R}^{\geq 0}$ is the quantifier function of the stability which depends only on $A$. Thus $Q(g)$ and $\rho$, $D-$fellow-travel each other where $D$ depends only on $A$ and the topological type of $S$. Moreover, since both $Q(g)$ and $\rho$ are quasi-geodesics with quantifiers depending only on the topological type of the surface $S$, there is a coarse parameter map $N$ from $[m,n]$ to $[a,b]$ such that 
\begin{equation}\label{eq : timecor}|N(i)-N(j)|\asymp_{K,C} |i-j|,\end{equation}
where the constants $K\geq 1$ and $C\geq 0$ depend only on $A$; see \cite[\S 5.3]{wpbehavior}.

The following theorem from \cite[\S 8]{wpbehavior} provides us a {\it WP geodesic ray with a prescribed itinerary} in the Teichm\"{u}ller space.

\begin{thm}\textnormal{(Infinite ray with prescribed itinerary)}\label{thm : inftyray}
Given $A,R,R'>0$, there are constants $\bar{w}=\bar{w}(A,R,R',\epsilon)$ and $\bar{\epsilon}=\bar{\epsilon}(A,R')$ with the following properties.
 
 Let $(\nu^{-},\nu^{+})$ be an $A-$narrow pair. Let $\rho$ be a hierarchy path between $\nu^{-}$ and $\nu^{+}$. Let $r_{\nu^{\pm}}:[0,\infty)\to\Teich(S)$ be the infinite WP geodesic ray as in Lemma \ref{lem : inftyraynu}. Suppose that $Z$ a large component domain of $\rho$ has $(R,R')-$bounded combinatorics over an interval $[m',n']\subset J_{Z}$ with $n'-m'>2\bar{w}$. Let $a'\in N(m'+\bar{w})$ and $b'\in N(n'-\bar{w})$. Then for every $t \in [a',b']$ we have
   \begin{enumerate}
 \item $\ell_{\gamma}(r_{\nu^{\pm}}(t))>\bar{\epsilon}$ for every $\gamma \notin \partial{Z}$, and
  \item $\ell_{\alpha}(r_{\nu^{\pm}}(t))\leq \epsilon $ for every $\alpha \in \partial{Z}$
 \end{enumerate}
 Moreover, if $Z_{1}$ and $Z_{2}$ are subsurfaces as above, $n'_{1}<m'_{2}$ implies that $b'_{1}<a'_{2}$.
 \end{thm}

\subsection{Bounding annular coefficients}\label{subsec : bdanncoeff}

Let $r_{\nu^{\pm}}:[0,\infty)\to \Teich(S)$ be a WP geodesic ray with prescribed itinerary where the end invariant $(\nu^{-},\nu^{+})$ is $A-$narrow and has $R'-$bounded annular coefficients. To goal of this section is to prove Lemma \ref{lem : rvisitingU}. Where show that over any long enough subinterval of $[0,\infty)$, there is a subinterval of definite length over which $r_{\nu^{\pm}}$ is in some region $U$ of the Teichm\"{u}ller space ($U$ regions were defined in $\S 2.2$). This combinatorial control will be used in $\S$\ref{sec : contractiongeodthin}.

First we recall two properties of hierarchy paths which will be used in this section. For an extended list of properties of hierarchy paths see \cite[\S 2]{wpbehavior}, \cite[\S 2]{bmm2}. 
\begin{thm}\label{thm : hr}
There exist positive constants $M_{2}$ and $M_{3}$ depending only on the topological type of $S$ with the following properties. Let $\rho:[m,n]\to P(S)$ be a hierarchy path between partial markings or laminations $\mu^{-}$ and $\mu^{+}$. Let $Z$ be a component domain of $\rho$ and let $J_{Z}=[j^{-},j^{+}]$. Then 
\begin{enumerate}
\item\label{hr : lfbddproj} $d_{Z}(\rho(i),\rho(j^{-}))\leq M_{2}$ for any $i\leq j^{-}$, and $d_{Z}(\rho(i),\rho(j^{+}))\leq M_{2}$ for any $i\geq j^{+}$.
\item (No backtracking) \label{hr : nobacktrack} Let $i,j,k\in [m,n]$ with $i\leq j\leq k$. Then for any subsurface $Y\subseteq S$, 
$$d_{Y}(\rho(i),\rho(k))+M_{3}\geq d_{Y}(\rho(i),\rho(j))+d_{Y}(\rho(j),\rho(k)).$$
\end{enumerate}
\end{thm}

\begin{lem} \label{lem : bddcombpart}
Given an increasing function $F: \mathbb{R}^{\geq0} \to \mathbb{R}^{\geq 0}$, there is an $L>0$, depending only on $F$ and the topological type of $S$ with the following property. Let $\rho: [m,n] \to P(S)$ be a hierarchy path. Suppose that a subinterval $[m',n']\subseteq [m,n]$ has the property that for any subsurface $Z$ and any $R>0$, if $Z$ has non-annular $R-$bounded combinatorics over a subinterval $[k,l]\subset [m',n']$, then $d_{Z}(\rho(l),\rho(k))\leq F(R)$. Then we have that $n'-m'\leq L$.
\end{lem}

\begin{proof}
Let $\xi(S)$ be the complexity of the surface $S$. For each $i=0,1,...,\xi(S)$ define the constant
\begin{equation}\label{eq : xi}x_{i}=\max\{d_{Y}(\rho(m'),\rho(n')):Y \subseteq S\;\text{ is non-annular and}\; \xi(Y)\leq i\}.\end{equation} 
Note that for $i=0$, the set on the right-hand side is empty and we define $x_{0}=0$.
 \begin{claim}\label{claim : dZm'n'bd}
 Let $i\in\{1,...,\xi(S)\}$. For any essential subsurface $Z\subseteq S$ with $\xi(Z)=i$ we have 
 $$d_{Z}(\rho(m'),\rho(n'))\leq F(x_{i-1})+2M_{2}+4.$$  
 \end{claim}
 The proof of the claim is by contradiction. Suppose that the claim does not hold. Then there is a non-annular subsurface $Z$ with $\xi(Z)=i$ such that 
 \begin{equation}\label{eq : dZm'n'>F}d_{Z}(\rho(m'),\rho(n'))>F(x_{i-1})+2M_{2}+4.\end{equation} 

 \begin{claim}\label{claim : JZint[m',n']}
Suppose that (\ref{eq : dZm'n'>F}) holds. Then $J_{Z}\cap[m',n']\neq\emptyset$.
 \end{claim}
 Let $J_{Z}=[j^{-},j^{+}]$. To get a contradiction, suppose that the claim does not hold. Then we have that $[j^{-},j^{+}]\cap [m',n']= \emptyset$. This implies that either $n'<j^{-}$ or $j^{+}<m'$. 
 
 First suppose that $n'< j^{-}$, then by Theorem \ref{thm : hr}(\ref{hr : lfbddproj}) we have $d_{Z}(\rho(n'),\rho(j^{-}))\leq M_{2}$. Moreover, since $m'<n'< j^{-}$, by Theorem \ref{thm : hr}(\ref{hr : lfbddproj}), we have $d_{Z}(\rho(m'),\rho(j^{-}))\leq M_{2}$. The last two inequalities combined with the triangle inequality give us 
 $$d_{Z}(\rho(m'),\rho(n'))\leq 2M_{2}+\diam_{Z}(\rho(j^{-}))\leq 2M_{2}+2.$$
where the second inequality above follows from Lemma \ref{lem : diamproj}. But this contradicts the bound (\ref{eq : dZm'n'>F}) we assumed to hold. 

Now suppose that $m'>j^{+}$, then by Theorem \ref{thm : hr}(\ref{hr : lfbddproj}) we have $d_{Z}(\rho(m'),\rho(j^{+}))\leq M_{2}$. Moreover, since $n'>m'>j^{+}$, by Theorem \ref{thm : hr}(\ref{hr : lfbddproj}) we have $d_{Z}(\rho(n'),\rho(j^{+}))\leq M_{2}$. The last two inequalities combined with the triangle inequality imply that 
$$d_{Z}(\rho(m'),\rho(n'))\leq 2M_{2}+2,$$
 which again contradict the bound (\ref{eq : dZm'n'>F}). This finishes the proof of Claim \ref{claim : JZint[m',n']}.
 \medskip
 
By Claim \ref{claim : JZint[m',n']} we have that $J_{Z}\cap[m',n']\neq\emptyset$. Define the interval 
$$[k,l]=J_{Z}\cap[m',n'].$$
 We proceed to show that 
 \begin{equation}\label{eq : dZkl>}d_{Z}(\rho(k),\rho(l))\geq F(x_{i-1}).\end{equation} 
 
 If $[m',n']\subseteq J_{Z}$, then $[k,l]=[m',n']$ and (\ref{eq : dZkl>}) follows immediately from the bound (\ref{eq : dZm'n'>F}). 
 
 If $[k,l]=[j^{-},n']$, then by the triangle inequality and Lemma \ref{lem : diamproj},
 \begin{eqnarray*}
 d_{Z}(\rho(m'),\rho(j^{-}))+d_{Z}(\rho(j^{-}),\rho(n'))&\geq& d_{Z}(\rho(m'),\rho(n'))-\diam_{Z}(\rho(j^{-}))\\
 &\geq& d_{Z}(\rho(m'),\rho(n'))-2.
 \end{eqnarray*}
Moreover $m'<j^{-}$ so by Theorem \ref{thm : hr}(\ref{hr : lfbddproj}) we have
$$d_{Z}(\rho(m'),\rho(j^{-}))\leq M_{2}.$$
 Now the above two inequalities and the inequality (\ref{eq : dZm'n'>F}) together imply that the inequality
$$d_{Z}(\rho(j^{-}),\rho(n'))\geq F(x_{i-1})+M_{2},$$ 
holds, from which the bound (\ref{eq : dZkl>}) follows.
 
 If $[k,l]=[m',j^{+}]$, then as above we can get the bound 
 $$d(\rho(m'),\rho(j^{+}))\geq F(x_{i-1})+M_{2},$$
  which gives us the bound (\ref{eq : dZkl>}). 
 
 Finally, if $J_{Z}\subseteq [m',n']$, then Theorem \ref{thm : hr}(\ref{hr : lfbddproj}) implies that 
 \begin{eqnarray*}
 d_{Z}(\rho(m'),\rho(j^{-}))&\leq& M_{2},\;\text{and}\\
  d_{Z}(\rho(n'),\rho(j^{+}))&\leq& M _{2}.
  \end{eqnarray*} 
  The above two inequalities and (\ref{eq : dZm'n'>F}) combined with the triangle imply that
   $$d_{Z}(\rho(j^{-}),\rho(j^{+}))\geq  F(x_{i-1}),$$
    which is the bound (\ref{eq : dZkl>}). This completes establishing of (\ref{eq : dZkl>}).
 \medskip
 
Now by the setup of the constants $x_{i}$ in (\ref{eq : xi}) the subsurface $Z$ has non-annular $x_{i-1}-$bounded combinatorics over the subinterval $[k,l]\subseteq[m',n']$ and by (\ref{eq : dZkl>}) we have
$$d_{Z}(\rho(k),\rho(l))\geq F(x_{i-1}).$$
 But this contradicts the assumption of the lemma. Claim \ref{claim : dZm'n'bd} follows from this contradiction.
\medskip
 
 By Claim \ref{claim : dZm'n'bd}, for each $i\in\{1,..., \xi(S)\}$, we have
 $$\max\{d_{Z}(\rho(m'),\rho(n')):\xi(Z)=i\}\leq F(x_{i-1})+2M_{2}+4,$$
  moreover from the setup of the constants $x_{i}$ in (\ref{eq : xi}), it is clear that 
$$x_{i}=\max\{x_{i-1},\max\{d_{Z}(\rho(m'),\rho(n')):\xi(Z)=i\}\},$$
 so we have the bound
\begin{equation}\label{eq : xiub}x_{i} \leq \max\{x_{i-1},F(x_{i-1})+2M_{2}+4\}.\end{equation}
To simplify the notation we define the function 
$$f(x)=\max\{x,F(x)+2M_{2}+4\}.$$
 By (\ref{eq : xiub}), $x_{i}\leq f(x_{i-1})$. Also we have that $x_{0}=0$. Then we may inductively show that
\begin{equation} \label{eq : indbd} x_{i} \leq f^{i}(0),\end{equation}
for each $i\in\{1,...,\xi(S)\}$, where $f^{i}$ denotes the $i-$th composition of $f$ with itself.

Since $f$ is an increasing function and $f(0)\geq 0$, we have that
$$\max\{f^{i}(0) : i=1,2,...,\xi(S)\}=f^{\xi(S)}(0).$$
 Then using the bound (\ref{eq : indbd}) we see that
$$\max\{x_{i}:i=1,2,...,\xi(S)\}\leq f^{\xi(S)}(0).$$
 So by the setup of the constants $x_{i}$ in (\ref{eq : xi}) we have 
\begin{equation}\label{eq : dYrm'rn'ub}d_{Y}(\rho(m'),\rho(n'))\leq f^{\xi(S)}(0)\end{equation}
 for all non-annular subsurfaces $Y\subseteq S$. Let $A_{1}=\max\{M_{1},f^{\xi(S)}(0)\}$ be the threshold constant in the distance formula (\ref{eq : dsf}). Let $C_{1}$ be the additive constant corresponding to $A_{1}$. Then by (\ref{eq : dYrm'rn'ub}) we have
$$d(\rho(n'),\rho(m')\leq C_{1}.$$
 Moreover, $\rho$ is a $(k,c)-$quasi-geodesic where $k$ and $c$ depend only on the topological type of $S$, so we obtain the upper bound $L=kC_{1}+kc$ for $n'-m'$. Note that $L$ only depends on $F$ and the topological type of $S$. The proof of the lemma is complete.
\end{proof}
Let $\bar{\epsilon}=\bar{\epsilon}(A,R')$ be the constant from Theorem \ref{thm : inftyray}.
\begin{lem}\label{lem : rvisitingU}
Given positive constants $A, d, R'$ and $\epsilon\leq\bar{\epsilon}$, there are constants $L'>0$ and $\bar{w}>0$ with the following properties. Let $(\nu^{-},\nu^{+})$ be an $A-$narrow pair with $R'-$bounded annular coefficients. Let $\rho:[m,n]\to P(S)$ be a hierarchy path between $\nu^{-}$ and $\nu^{+}$. Let $r_{\nu^{\pm}}:[a,b) \to \Teich(S)$ be a WP geodesic ray with prescribed itinerary and end invariant $(\nu^{-},\nu^{+})$. Then for any subinterval $[m',n']\subseteq [m,n]$ with $m'-n' \geq L'$, there are a subinterval $[k,l]\subset [m',n']$ and a large component domain $Z$ of $\rho$ such that for $t^{-}\in N(k+\bar{w})$ and $t^{+}\in N(l-\bar{w})$ ($N$ is the parameter correspondence map) we have $r_{\nu^{\pm}}(t)\in U_{\epsilon,\bar{\epsilon}}(\partial{Z})$ for every $t\in [t^{-},t^{+}]$. Moreover, $t^{+}-t^{-}\geq d$. 
\end{lem}
\begin{proof}
 Fix a threshold constant $A_{1}\geq M_{1}$ for the distance formula (\ref{eq : dsf}) and let $K_{1},C_{1}$ be the corresponding constants. Note that the hierarchy path $\rho$ is a $(k,c)-$quasi-geodesic where $k,c$ depend only on the topological type of the surface. Let $K_{2}=K_{1}k$ and $C_{2}=K_{1}k(C_{1}+c)$. Let $K,C$ be the constants for $N$ from (\ref{eq : timecor}). Let $M_{3}$ be the constant from Theorem \ref{thm : hr}(\ref{hr : nobacktrack}). Let $\bar{w}$ be the constant from Theorem \ref{thm : inftyray}. Define the function
$$F(x)=\max\{K_{2}\big(2\bar{w}(A,x,R'+2M_{3},\epsilon)+(Kd+KC)\big)+C_{2},A\}.$$
Now let $L'$ be the constant from Lemma \ref{lem : bddcombpart} for the the function $F$ defined above. For any subinterval $[m',n']\subseteq [m,n]$ with $n'-m'\geq L'$ by the contrapositive of Lemma \ref{lem : bddcombpart}, there are $R>0$, a subsurface $Z$ and an interval $[k,l] \subseteq [m',n']\cap J_{Z}$ such that 
$$d_{Z}(\rho(k),\rho(l))> F(R),$$
 and $Z$ has non-annular $R-$bounded combinatorics over $[k,l]$. Since $\rho$ is $A-$narrow and $F(R)\geq A$, the subsurface $Z$ is a large subsurface. Thus for any non-annular subsurface $Y$, either $Y\subseteq Z$ or $Y\pitchfork Z$ holds. If $Y\pitchfork Z$, then since $\partial{Z}\subset\rho(l),\rho(k)$, we have $d_{Y}(\rho(k),\rho(l))\leq 2$. If $Y\subseteq Z$, then by the non-annular $R-$bounded combinatorics we have $d_{Y}(\rho(k),\rho(l))\leq R$. Then by the distance formula (\ref{eq : dsf}) and the fact that $Z$ has $R-$bounded combinatorics over $[k,l]$, we have
$$d(\rho(l),\rho(k))\geq \frac{1}{K_{1}}\{d_{Z}(\rho(l),\rho(k))\}_{A_{1}}-C_{1}.$$
As we saw above $d_{Z}(\rho(k),\rho(l))> F(R)$, so 
 $$d(\rho(l),\rho(k))\geq k(2\bar{w}(A,R,R'+2M_{3},\epsilon)+(Kd+KC))+kc.$$
  Moreover, $\rho$ is a $(k,c)-$quasi-geodesic, so 
  $$k-l>2\bar{w}(A,R,R'+2M_{3},\epsilon)+Kd+KC.$$
 Furthermore, by the assumption of the lemma, $d_{\gamma}(\nu^{-},\nu^{+})\leq R'$
 for any $\gamma \in \mathcal{C}_{0}(S)$. Then by Theorem \ref{thm : hr}(\ref{hr : nobacktrack}), 
 $$d_{\gamma}(\rho(k),\rho(l))\leq R'+2M_{3}.$$
  Therefore, $Z$ has $(R,R'+2M_{3})-$bounded combinatorics over $[k,l]$. Set $\bar{w}=\bar{w}(A,R,R'+2M_{3},\epsilon)$. Then the lemma follows from applying Theorem \ref{thm : inftyray} to the interval $[k,l]$. 
\end{proof}

\section{Variation of geodesics}\label{sec : ruledrec}

Let $\mathcal{X}$ be a geodesically convex, negatively curved, Riemannian manifold; for example the Teichm\"{u}ller space equipped with the WP metric. Let $g:[a,b]\to \mathcal{X}$ be a geodesic segment, and let $\pi:\mathcal{X}\to g$ be the nearest point projection from $\mathcal{X}$ to $g$. In the following proposition we collect some important facts about the map $\pi$ which we need:
\begin{prop}\label{prop : proj map}
Suppose that $p$ is a point in $\mathcal{X}$ so that $\pi(p)$ is in the interior of $g([a,b])$ or is an end point $g(a)$ or $g(b)$ which is the nearest point to $p$ on a slightly longer geodesic segment containing $g([a,b])$. Let $\zeta$ be a geodesic segment connecting $p$ and $\pi(p)$. We have
\begin{enumerate}
\item  The projection map $\pi$ is continuous at $p$. 
\item If $p$ is not on $g([a,b])$, then $\pi$ is smooth at $p$.
\item The distance function $d(p,\pi(p))$ is smooth.
\end{enumerate}
\end{prop}
\begin{proof}
Part (1) is \cite[Lemma 3.2]{mfldnc}. Part (2) follows from \cite[Proposition 1.7]{chebcomparison}. Part (3) follows from the lemma on the first page of \cite{regdist}. Part (4) is \cite[Theorem 1]{regdist}. 
\end{proof}

Let $g_{1}:[a,b]\to \mathcal{X}$ and $g_{2}:[c,d]\to \mathcal{X}$ be two geodesic segments parametrized by arclength. Let $g_{1}(t)$ be a point whose nearest point projection on $g_{2}([c,d])$ is in the interior of $g_{2}([c,d])$. Note that $g_{1}$ is a smooth path, moreover by Proposition \ref{prop : proj map}(4) the end point of $[g_{1}(t),\pi(g_{1}(t))]$ varies smoothly. Therefore $[g_{1}(t),\pi(g_{1}(t))]$ is a smooth family of geodesic segments with respect to $t\in[a,b]$. Moreover, since $\mathcal{X}$ is a negatively curved manifold, any two of the geodesic segments $[g_{1}(t),\pi(g_{1}(t))]$ fellow travel all the way. Finally, let $f:[a,b]\to g_{2}([c,d])$ be the reparametrization of $g_{2}$ which maps $t\in[a,b]$ to $\pi(g_{1}(t))$. Then $f$ is smooth.

 Let 
 $$\triangle=\{(t,s): t\in[a,b]\;\text{and}\; s\in [0,\lambda(t)]\},$$
  where $\lambda(t)$ is the length of the geodesic segment $[g_{1}(t),\pi(g_{1}(t))]$. For any $t\in[a,b]$, let $u(t,s)$ be the arclength parametrization of $[g_{1}(t),g_{1}(\pi(t))]$ where $s$ is the parameter and $u(t,0)=\pi(g_{1}(t))$. Then since $[g_{1}(t),g_{1}(\pi(t))]$ is a smooth family of geodesics, the map $(t,s)\mapsto u(t,s)$ defines a smooth map $u:\triangle\to\mathcal{X}$. For any $t\in [a,b]$, 
  $$\zeta(t):=u(t,s)$$
   is a geodesic parametrized by $s$. Moreover, by Proposition \ref{prop : proj map}(3), for any $s>0$, $u(t,s)$ is a smooth path parametrized by $t$. Then $u$ is a variation of geodesics $[g_{1}(t),\pi(g_{1}(t))]$ and $\frac{\partial u}{\partial t}$ is a Jacobi field. 
 
\begin{figure}
\centering
\scalebox{0.2}{\includegraphics{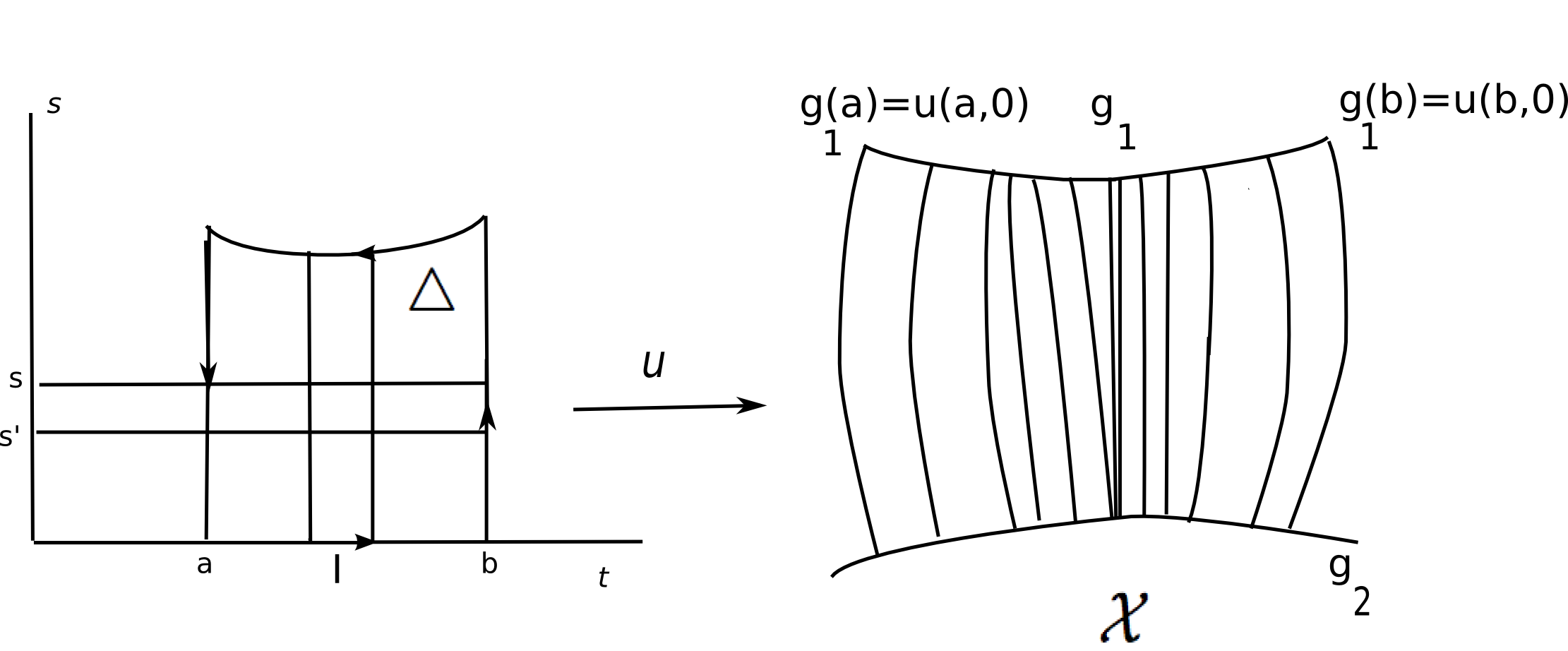}}
\caption{The variation of geodesics $u$ is the family of geodesic segments connecting points of the geodesic segment $g_{1}$ to their nearest points on the geodesic segment $g_{2}$. We pull back of the metric of $\mathcal{X}$ via $u$ to the region $\triangle$ in the $(t,s)$ plane.}
\label{fig : negcurv}
\end{figure}

 For any $s$, $(.,s)$ defines a vertical coordinate line and for any $t$, $(t,.)$ defines a horizontal coordinate line in $\triangle\subset\mathbb{R}^{2}$; see Figure \ref{fig : negcurv}.

We pull back the metric of $\mathcal{X}$ to $\triangle$ via $u$. Let the interval $I\subset [a,b]$ and $s,s'>0$ with $s'<s$ be such that $I\times s \subset \triangle$ and $I\times s' \subset \triangle$. The main result of this section is Lemma \ref{lem : lengthcurv}, where we prove that a difference between the length of $I\times s$ and $I\times s'$ gives rise to a definite total Gaussian curvature of the region $I\times [s',s]\subset \triangle$. For this purpose, using the Gauss-Bonnet formula, we prove the formula (\ref{eq : gbsquare}).
 
 Let $k_{g}(t,s)$ be the geodesic curvature of the path $u(.,s)$ at $u(t,s)$, where the normal vector $\vec{n}$ of the path is the unit vector orthogonal to $u$ which satisfies $\langle \vec{n}(u(t,s)),\frac{\partial{u}}{\partial{s}}(t,s)\rangle <0$. 
 
\begin{thm}\label{thm : pullbackm}
The pull back metric has the following properties:
\begin{enumerate}[(I)]
\item \label{pmetric : vergeod} Each vertical coordinate line is a geodesic.
\item \label{pmetric : verhor} Vertical and horizontal coordinate lines intersect each other orthogonally.
\item \label{pmetric : k} There is a $k<0$ such that the Gaussian curvature $\kappa$ is $\leq k$.
\item \label{pmetric : distinc} The distance between any two vertical coordinate lines is increasing in $s$, so for any $I\subseteq [a,b]$ if $s\geq s'$ then $\length(I\times s)\geq \length(I\times s')$.
\item \label{pmetric : kg} The function $-k_{g}$ is non-negative along each horizontal coordinate line.
\end{enumerate}
\end{thm}
\begin{proof}
 Part (\ref{pmetric : vergeod}) follows from the fact that the path $u(t,.)$ is a geodesic in $\mathcal{X}$. Moreover, the paths $u(t,.)$ and $u(.,s)$ intersect orthogonally in $\mathcal{X}$, so part (\ref{pmetric : verhor}) follows. 

The family of geodesic segments $[g_{1}(t),\pi(g_{1}(t))]$, where $t\in[a,b]$, is a compact subset of $\mathcal{X}$. Thus the sectional curvatures of the metric of $\mathcal{X}$ at any point of $u(\Delta)$ are bounded above by some $k <0$. Moreover exercise II.16 on page 104 of \cite{chaveldg} asserts that: The Gaussian curvature of the ruled surface obtained from a geodesic variation at any point is less than or equal to the sectional curvature of the tangent plane of the surface at that point. Therefore part (\ref{pmetric : k}) follows. 

As we saw above the sectional curvature of $\mathcal{X}$ at any point of $u(\Delta)$ is bounded above by $k<0$. It is a standard fact that for $u(t,s)$ a smooth geodesic variation as above in a manifold with constant negative sectional curvature $k$ the length of a Jacobi field defined by $\frac{\partial u}{\partial t}$ is increasing along $u(t,s)$ as $s$ increases; In fact this follows from the explicit formula for the Jacobi fields in a constant negative curvature manifold given in \cite[\S II.5]{chaveldg}. Then the Rauch comparison theorem for the region in $\mathcal{X}$ with sectional curvatures bounded above by $k$ and a manifold with constant curvature $k$, implies that the length of the Jacobi field $\frac{\partial u}{\partial t}$ is increasing. For more detail about the Rauch comparison theorem see e.g. $\S 11$ of Chapter 1 of \cite{chebcomparison}. So we may conclude that 
$$\frac{\partial u}{\partial t}|_{(t,s)}\geq \frac{\partial u}{\partial t}|_{(t,s')}$$
 for any $t$ and any $s,s'$ with $s\geq s'$. Moreover, for any $s$,
 $$\length(u(I\times s))=\int_{I}\frac{\partial u}{\partial t}|_{(t,s)} dt.$$ Then using the inequality above, we have
$$\length(u(I\times s))\geq\length(u(I\times s')),$$
which is part (\ref{pmetric : distinc}). 

Let $I= [a',b']$ be a subinterval of the interval $[a,b]$. We may consider $u$ as a variation of paths $u|_{I\times s}$. Then by the first variation of arc-length formula (see e.g. page 4 of \cite{chebcomparison}) we have
$$\frac{d}{ds}\length(I\times s)= \langle \frac{\partial{u}}{\partial{s}}, \frac{\partial{u}}{\partial{t}} \rangle\Big|_{a'}^{b'}-\int_{I\times s} \langle \frac{\partial{u}}{\partial{s}}, \nabla_{\frac{\partial{u}}{\partial{t}}}\frac{\partial{u}}{\partial{t}} \rangle dt.$$
 By (\ref{pmetric : verhor}) we have that $\langle \frac{\partial{u}}{\partial{t}},\frac{\partial{u}}{\partial{s}}\rangle \equiv 0$, so the first term on the right hand-side of the above formula is $0$. Moreover by definition $k_{g}:=\langle \frac{\partial{u}}{\partial{s}}, \nabla_{\frac{\partial{u}}{\partial{t}}}\frac{\partial{u}}{\partial{t}}\rangle$. Thus we have 
\begin{equation}\label{eq : fvarclen}\frac{d}{ds}\length(u(I\times s))=\int_{I\times s}-k_{g}(t,s)dt.\end{equation} 
Moreover, by (\ref{pmetric : distinc}), 
$$\frac{d}{ds}\length(u(I\times s))\geq 0.$$
 Therefore 
$$\int_{I\times s}-k_{g}\geq 0.$$
 Since $I$ was an arbitrary subinterval of $[a,b]$, from the fact that the integral of $-k_{g}$ is positive over $I$, we can conclude that the continuous function $-k_{g}$ is a non-negative function on the interval $[a,b]$. This is part (\ref{pmetric : kg}). 
\end{proof}

\subsection{Gauss-Bonnet formula} For any $s>0$, by Theorem \ref{thm : pullbackm}(\ref{pmetric : kg}), the integral of the function $-k_{g}$ with respect to $dt$ along the horizontal coordinate line $(.,s)$ defines a positive measure $m_{s}$ on $(.,s)$. 

Suppose that $V \times W \subset \triangle$. Let $I\subseteq V$ and $[s',s]\subset W$ be subintervals. By the Gauss-Bonnet formula (page 242 of \cite{chaveldg}) for the pull back metric on $I\times [s',s]$ we have
  $$\int\int_{I\times[s',s]} \kappa dA+\int_{\partial({I\times[s',s]})} k_{g}dt=2\pi-\sum_{i}\theta_{i},$$
 where the orientations of the boundary curves in the line integral above are shown in Figure \ref{fig : negcurv}. By Theorem \ref{thm : pullbackm}(\ref{pmetric : vergeod}) each $t\times (s',s)$ is a geodesic. The sum above is taken over the exterior angles at the four corners of  $I\times[s',s]$. Moreover, by Theorem \ref{thm : pullbackm}(\ref{pmetric : verhor}) each exterior angle is equal to $\frac{\pi}{2}$, so the exterior angles add up to $2\pi$. Therefore we get 
\begin{equation}\label{eq : gb'}-\int_{I\times s}k_{g}dt+\int_{I\times s'}k_{g}dt+\int\int_{I\times [s',s]} \kappa dA=0.\end{equation}
After rearrangement of terms above we have
\begin{equation}\label{eq : kgcurv}\int_{I\times s'} k_{g}dt=\int_{I\times s}k_{g}dt+ \int\int_{I\times [s',s]}-\kappa dA.\end{equation}

\begin{claim}\label{claim : kglim}
The limit of the left hand-side of (\ref{eq : kgcurv}) exists as $s'\to 0$. 
\end{claim}
First, by Theorem \ref{thm : pullbackm}(\ref{pmetric : kg}) for every $s'>0$, $k_{g}\leq 0$, so the left hand-side of (\ref{eq : kgcurv}) is non-positive. Thus the right-hand side of (\ref{eq : kgcurv}) is non-positive as well. Second, the first integral on the right-hand side of (\ref{eq : kgcurv}) does not depend on $s'$. Moreover, by Theorem \ref{thm : pullbackm}(\ref{pmetric : k}), $\kappa\leq0$, so the second integral on the right-hand side is increasing as $s'\to 0$. Thus the right-hand side is increasing as $s'\to 0$. As a result the left-hand side is a non-positive function of $s'$ which increases as $s'\to 0$. Thus the limit of the left-hand side exits as $s'\to 0$. Furthermore, since the equality holds for every $s'>0$, it holds at $s'=0$ as well. The claim is proved.
\medskip

Since Claim \ref{claim : kglim} holds for every subinterval of $I$, the weak$^{*}$ limit of the measures $m_{s}$ exists. Denote by $m$ the measure obtained as the weak$^{*}$ limit of the measures $m_{s}$. Then taking limit as $s'\to 0$ in (\ref{eq : kgcurv}) we obtain the formula 
\begin{equation}\label{eq : gbsquare}\int\int_{I\times [0,s]}\kappa dA-\int_{I\times 0} m-\int_{I\times s} k_{g}dt=0.\end{equation}

Consider the regions in $\triangle$ below the horizontal lines; see Figure \ref{fig : negcurv}. Define the measures $m_{s}$ on the horizontal coordinate lines $(.,s)\cap\triangle$ by integrating $-k_{g}$. Let the measure $m$ be the limit of the measures $m_{s}$ as $s\to 0$. Then the argument given above to prove the formula (\ref{eq : gbsquare}) gives us the following formula
\begin{equation} \label{eq : gb}\int\int_{\triangle}\kappa dA-\int_{[a,b]\times 0}m=2\pi-\sum_{i}\theta_{i},\end{equation}
where the sum is taken over the exterior angles of the boundary of the region $\triangle$.

\subsection{A length versus total curvature estimate} 

As before, let $V\times W\subset \triangle$, and also let $I\subseteq V$ and $[s,s']\subset W$ be subinterval. Then define
  \begin{equation}\label{eq : phi}\phi(s,s') = \Big\lbrace 
  \begin{array}{cl}
  \int\int_{I\times [s',s]} -\kappa dA -\int_{I\times s'}k_{g}dt  & \text{if}\; s' >0\\
  \int\int_{I\times [0,s]} -\kappa dA +\int_{I\times 0}m & \text{if}\; s' = 0.\\
  \end{array}\end{equation}
 For any $s\in W$ let 
 $$l_{s}=\length(u(I\times s)).$$
 Then $l_{s}$ is positive and by Theorem \ref{thm : pullbackm}(\ref{pmetric : distinc}) decreases as $s\to 0$, therefore the limit $\lim_{s\to 0}l_{s}$ exists. We denote $l_{0}=\lim_{s\to 0} l_{s}$. 
 \begin{remark}
 Note that $l_{0}$ is not necessarily equal to $\length(u(I\times 0))$. 
 \end{remark}
 \begin{lem} \label{lem : lengthcurv} For any $s,s'\in W$ with $s'\leq s$, we have
  $$\phi(s,s') \geq \frac{l_{s}-l_{s'}}{s-s'}.$$
   \end{lem}
 \begin{proof}
First  we prove the following integral formula:
\begin{eqnarray}\label{eq : integralformula}
l_{s}-l_{s'}&=&\int_{s'}^{s} \frac{dl_{r}}{dr}  dr\\ \nonumber
&=&\int_{s'}^{s} \Big(\int_{I\times r} -k_{g}dt\Big)dr\\ \nonumber
&=&\int_{s'}^{s} \Big(\int\int_{I\times [s',r]} -\kappa dA + \int_{I\times s'}-k_{g}dt\Big)dr.\nonumber
\end{eqnarray}
The first equality is the Fundamental theorem of calculus applied to $l_{r}$ as a function of $r$. The second equality holds by the first variation of arc-length formula (\ref{eq : fvarclen}) for the orthogonal variation of paths $u|_{I\times s}$. The last equality holds by the Gauss-Bonnet formula (\ref{eq : gb'}). 
 
Now suppose that $s>0$ and $r,s'\in(0,s]$ such that $r\geq s'$. By the definition of $\phi$ in (\ref{eq : phi}), $\phi(s,r)\leq \phi(s,s')$. So we have 
\begin{center}$\int_{s'}^{s}\phi(s,r)dr\leq \int_{s'}^{s}\phi(s,s')dr=\phi(s,s')(s-s'),$ \end{center}
then using (\ref{eq : integralformula}) we have
$$l_{s}-l_{s'} \leq \phi(s,s')(s-s').$$
This proves the lemma for any $s'>0$. 

 We proceed to show that the lemma holds at $s'=0$ as well. For any fixed $r<s$ by (\ref{eq : kgcurv}) the integrand in the last line of (\ref{eq : integralformula}) does not depend on $s'$. So the limit as $s'\to 0$ of the last line of (\ref{eq : integralformula}) is $\int_{0}^{s}\phi(0,r)dr$. Thus $l_{s}-l_{0}=\int_{0}^{s}\phi(0,r)dr$. Moreover, by (\ref{eq : phi}), for $s\geq r$, $\phi(s,0)\geq \phi(r,0)$. So we have
$$l_{s}-l_{0}\leq\int_{0}^{s}\phi(r,0)dr \leq s\phi(s,0).$$
Therefore, the lemma holds at $s'=0$ as well. The proof of the lemma is complete.
\end{proof}

\begin{remark}
If the variation of geodesics $u$ extends to a variation that contains $g_{1}([a,b])$ in its interior, then the second term of $\phi(0,s)$ in (\ref{eq : phi}) vanishes. The definition of the measure $m$ and the rather long discussion of this section meant to handle the possibility that the variation does not extend, for example when the nearest point projection onto $g_{1}([a,b])$ maps a subinterval to a point in the interior of $g_{1}([a,b])$.
\end{remark}

\section{Contraction property of WP geodesic segments}\label{sec : contractiongeodthin}
 Let $\mathcal{X}$ be a geodesically convex, negatively curved Riemannian manifold; For example the Teichm\"{u}ller space equipped with the WP metric. Let $g:[0,T]\to \mathcal{X}$ be a geodesic segment. Let $\pi:\mathcal{X}\to g([0,T])$ denote the nearest point projection map to $g([0,T])$.  At a point $p\in\mathcal{X}$ with $\pi(p)$ in the interior of $g([0,T])$, the map $\pi$ is smooth and has a linear derivative $d\pi:T_{p}\mathcal{X}\to T_{\pi(p))}g$ (by Proposition \ref{prop : proj map}(3)). At a point $p\in\mathcal{X}$, with $\pi(p)$ equal to either $g(0)$ or $g(T)$, only the directional derivatives of $\pi$ are defined. We denote the (directional) derivative of $\pi$ at $p$ by $d\pi_{p}:T_{p}\mathcal{X}\to T_{p}g$. The main result of this section is the following uniform (strict) contraction property of WP geodesic segments in certain regions of the Teichm\"{u}ller space.
\begin{thm}\label{thm : contnpt}
Given $\epsilon>0$ sufficiently small, $\bar{\epsilon},T$ and $b$ positive, there is a $\delta\in [0,1)$ with the following property. Let $\sigma$ be a possibly empty multi-curve such that the subsurface $S \backslash \sigma$ is a large subsurface. Let $g:[0,T]\to U_{\epsilon,\bar{\epsilon}}(\sigma)$ be a geodesic segment and let $\mathcal{N}_{b}(g([0,T]))$ be the $b-$neighborhood of $g([0,T])$. Then for every $p\in\Teich(S)\backslash\mathcal{N}_{b}(g([0,T]))$ and $v\in T_{p}\Teich(S)$ the inequality 
$$\frac{||d\pi_{p}(v)||^{2}}{||v||^{2}}\leq \delta$$
holds.  
\end{thm}

\begin{remark}
Compare this theorem with the contraction property of Teichm\"{u}ller geodesics which project to the thick part of the moduli space (co-bounded geodesics) proved by Minsky in \cite{contteichgeod}. Minsky uses the explicit description of flat surfaces along a Teichm\"{u}ller geodesic. But we use various estimates on the WP metric and its derivatives and a standard Jacobi field argument. 
\end{remark}

We start by collecting some of Wolpert's estimates for the Weil-Petersson metric and WP Levi-Civita covariant derivatives in the thin part of the Teichm\"{u}ller space. 

On a Riemannian manifold the Levi-Civita covariant derivative $\nabla$ (see e.g. \cite[\S I.5]{chaveldg}) is the unique covariant derivative which is
\begin{itemize}
\item   compatible with the Riemannian metric, i.e. for any smooth path $\zeta(t)$ and vector fields $V$ and $W$ along $\zeta$, $\frac{d}{dt}\langle V,W\rangle = \langle \nabla_{\dot{\zeta}}V,W \rangle+ \langle V,\nabla_{\dot{\zeta}}W \rangle$, and 
\item  torsion free, i.e. for any two vector fields $V$ and $W$, $\nabla_{V}W-\nabla_{W}V=[V,W]$. Here $[.,.]$ denotes the Lie bracket of vector fields.
\end{itemize}

Given a multi-curve $\sigma$ and $c_{0}>0$, let $\{\lambda_{\alpha},J\lambda_{\alpha},\grad \ell_{\beta}\}_{\alpha\in \sigma, \beta\in\varsigma}$ be the \textit{short and relative length frame field} in the region 
$$U_{c_{0}}(\sigma)=\{x\in \Teich(S) : \ell_{\alpha}(x)\leq c_{0}\; \text{for every}\; \alpha \in \sigma \}$$
 introduced by Wolpert in \cite{wolb}. Here $\varsigma$ is a marking on the surface $S\backslash\sigma$. $\lambda_{\alpha}=\grad \ell_{\alpha}^{1/2}$ and $\grad \ell_{\beta}$ are vector fields, and $J$ is the almost complex structure of the Teichm\"{u}ller space. 
 
 Recall  the completion stratum $\mathcal{S}(\sigma)$ corresponding to pinching the curves in the multi-curve $\sigma$ from $\S$\ref{subsec : wp}.

\begin{prop}\label{prop : wolpertmetric} We have the following estimates:
\begin{enumerate}
\item $\langle \lambda_{\alpha},\lambda_{\alpha'} \rangle=O(\ell_{\alpha}^{3/2}\ell_{\alpha'}^{3/2})$ for every $\alpha, \alpha'\in \sigma$ such that $\alpha\neq \alpha'$, and $ ||\lambda_{\alpha}||^{2}=\Theta(1)$ for every $\alpha \in \sigma$. 
\item $\langle \lambda_{\alpha},J\lambda_{\alpha'} \rangle=0$ for every $\alpha,\alpha'\in \sigma$.  
\item $\langle \grad \ell_{\beta},\lambda_{\alpha} \rangle=O(\ell_{\alpha}^{3/2})$ for every $\alpha \in \sigma$ and $\beta \in \varsigma$. 
\item $\langle \grad \ell_{\beta},J\lambda_{\alpha} \rangle=0$ for every $\alpha\in \sigma$ and $\beta\in \varsigma$.
\end{enumerate}
Moreover, for every $\beta,\beta'\in\varsigma$, $\langle \grad \ell_{\beta}, \grad \ell_{\beta'}\rangle$ is continuous in a neighborhood of the $\sigma-$stratum $\mathcal{S}(\sigma)\subset \overline{\Teich(S)}$. Here the constant of the $O$ notation and the constants of the $\Theta$ notation are uniform for $\ell_{\alpha}\leq c_{0}$.
\end{prop}
The estimates of the above proposition are established in Lemmas 3.12 and 4.2 of \cite{wolb}. See also Theorem 4.3 and Corollaries 4.3 and 4.4 of  \cite{wolb} where Wolpert puts these estimates together to get expansions for the WP metric near completion strata. 

 Note that the WP metric is the real part of a Hermitian metric on the Teichm\"{u}ller space. Therefore, $\langle V, W \rangle=\langle JV, JW \rangle$, where $J$ is the almost complex structure of the Teichm\"{u}ller space with the property that $J^{2}=-I$. Thus estimates for the remaining pairings of the vector fields of the short and relative length frame field follow from the ones listed in the above proposition.

\begin{prop}\label{prop : cdlambda}(\cite[Theorem 3.4]{wolextconn})
We have 
$$\nabla_{V}\lambda_{\alpha}=3\ell_{\alpha}^{-1/2}\langle J\lambda_{\alpha},V \rangle J\lambda_{\alpha}+O(\ell_{\alpha}^{3/2})$$
 where the constant of the $O$ notation is uniform for $\ell_{\alpha}\leq c_{0}$.
\end{prop}

\begin{prop}\label{prop : wolpertcovderiv} We have 
\begin{center}$\nabla_{\lambda_{\alpha}}\grad \ell_{\beta}=O(\ell_{\alpha}^{1/2})$ and $ \nabla_{J\lambda_{\alpha}}\grad\ell_{\beta}=O(\ell_{\alpha}^{1/2})$.\end{center} 
Moreover $ \nabla_{\grad\ell_{\beta'}}\grad \ell_{\beta}$ is continuous in a neighborhood of $\mathcal{S}(\sigma)\subset \overline{\Teich(S)}$. Here the constant of the $O$ notation is uniform for $\ell_{\alpha}\leq c_{0}$.
\end{prop}
 The estimates of this proposition are from Proposition 4.6 of \cite{wolextconn}. We also need the following estimates for the WP sectional curvatures. Let $V,U$ be two vector fields. Recall that the sectional curvature is defined by 
 $$\kappa(U,V)=\frac{\langle R(U,V)V,U \rangle}{|U\wedge V|^{2}},$$
  where $|U\wedge V|^{2}=||U||^{2}||V||^{2}-\langle U,V \rangle^{2}$.

\begin{prop}\label{prop : secstr}(\cite[Theorem 21]{wpcurv}) 
Let $\sigma$ be a multi-curve and $P$ be a pants decomposition, so that  $\sigma\subseteq P$. The diagonal curvature evaluations for $\alpha\in \sigma$ satisfy 
$$\langle R(\lambda_{\alpha},J\lambda_{\alpha})J\lambda_{\alpha}, \lambda_{\alpha}\rangle=3(16\pi\ell_{\alpha}^{3})^{-1}+O(\ell_{\alpha})$$
 and all of the remaining curvature evaluations are continuous in a neighborhood of $\mathcal{S}(\sigma)\subset \overline{\Teich(S)}$.
\end{prop}
In this proposition evaluations of the Riemann curvature tensor are in the frame $\{\lambda_{\alpha},J\lambda_{\alpha}\}_{\alpha \in P}$ (not the short and relative length frame). Moreover, Wolpert uses the convention that on $\mathcal{S}(\sigma)$ the evaluations of the Riemann curvature tensor $\langle R((J)\lambda_{\alpha},(J)\lambda_{\beta})(J)\lambda_{\gamma},(J)\lambda_{\delta}\rangle$ vanish in the situations:
\begin{itemize}
\item $\alpha \in \sigma$ and at least one of $\beta,\gamma$ and $\delta$ is distinct from $\alpha$,
\item  $\alpha,\beta,\gamma, \delta \in P-\sigma$ and not all of them on the same component of a Riemann surface with nodes represented in $\mathcal{S}(\sigma)$.
\end{itemize}

Define the bundles 
\begin{center}$N_{\sigma}=\spn\{\lambda_{\alpha}, J\lambda_{\alpha}\}_{\alpha\in \sigma}$ and $P_{\sigma}=\spn\{\grad_{\beta}\}_{\beta\in \varsigma}$\end{center}
 over $U_{c_{0}}(\sigma)$. Any vector field $V$ on the Teichm\"{u}ller space has a decomposition as $V=V_{N}+V_{P}$, where $V_{N}$ is a section of  $N_{\sigma}$ and $V_{P}$ is a section of $P_{\sigma}$. 

We proceed by proving two lemmas which elaborate the asymptotic product form of the WP metric.

\begin{lem}\label{lem : asympprod}
Let $\epsilon\leq c_{0}$. Let $p\in U_{\epsilon}(\sigma)$ and $v,w\in T_{p}\Teich(S)$ be two vectors with $||v||^{2}\leq 1$ and $||w||^{2}\leq 1$. Let 
\begin{center}$v=\sum_{\alpha\in \sigma}a_{\alpha}\lambda_{\alpha}+b_{\alpha}J\lambda_{\alpha}+\sum_{\beta\in \varsigma}c_{\beta}\grad\ell_{\beta}.$\end{center}
 Then 
\begin{enumerate}[(i)]
\item \label{asympprod : ao(1)} For any $\alpha \in \sigma$ and $\beta\in \varsigma$ the coefficients $a_{\alpha},b_{\alpha}$ and $c_{\beta}$ are $O(1)$. 
\item \label{asympprod : aprod} Suppose that $v\in P_{\sigma}$ and $w\in N_{\sigma}$. Then $\langle v,w\rangle|_{p} \to 0 $ as $\epsilon \to 0$.
\item \label{asympprod : anorm} The inequalities $||v||^{2}\geq \frac{1}{2}||v_{N}||^{2}$ and $||v||^{2}\geq\frac{1}{2}||v_{P}||^{2}$ hold.
\end{enumerate}
\end{lem}
\begin{proof}
Denote by $\overline{U_{\epsilon}(\sigma)}$ the closure of $U_{\epsilon}(\sigma)$ in $\overline{\Teich(S)}$. The short and relative length frame field $\{\lambda_{\alpha},J\lambda_{\alpha},\grad\ell_{\beta}\}_{\alpha\in \sigma, \beta\in \varsigma}$ on the region $U_{\epsilon}(\sigma)$ extends to a frame field on $\overline{U_{\epsilon}(\sigma)}$. Moreover the frame field is invariant under the action of the mapping class group, so projects to a frame field 
$$\mathbb{F}:=\{E_{i},E'_{i},F_{j}\}_{i=1,...,|\sigma|, j=1,...,|\varsigma|}$$
 on $\overline{\mathcal{M}(S)}$.  Since the WP metric is also invariant under the action of the mapping class group, the inner product of any two vectors in the frame is preserved. Thus, by the estimates of  Proposition \ref{prop : wolpertmetric} we have: 
 \begin{center}
 $||E_{i}||^{2}=\Theta(1)$, $||E'_{i}||^{2}=\Theta(1)$, $\langle E_{i},E'_{i'}\rangle=0$ if $i\neq i'$, \\
 $\langle E_{i},F_{j}\rangle=O(\epsilon^{3/2})$, $\langle F_{j},E'_{i}\rangle=0$ and $\langle F_{j},F_{j'} \rangle=\Theta(1)$.  
 \end{center}
 Note that the constant of the $O$ notation and the constants of the $\Theta$ notation depend only on $c_{0}$.

Denote by $\hat{p}$ the projection of the point $p$ to $\mathcal{M}(S)$. For any frame at $T_{\hat{p}}\mathcal{M}(S)$, the square of the lengths of the vectors of the frame and the inner product of any two of them together determine a point in $\mathbb{R}^{2|\sigma|+|\varsigma|}\times\mathbb{R}^{2|\sigma|+|\varsigma| \choose 2}$, and vice versa. So the set of frames at $\hat{p}$ are parametrized by $\mathbb{R}^{2|\sigma|+|\varsigma|}\times\mathbb{R}^{2|\sigma|+|\varsigma| \choose 2}$. The bounds above determine a compact subset $K(\hat{p})\subset \mathbb{R}^{2|\sigma|+|\varsigma|}\times\mathbb{R}^{2|\sigma|+|\varsigma| \choose 2}$. Moreover the frame $\mathbb{F}$ at $\hat{p}$ is in $K(\hat{p})$. Denote the projection of $\overline{U_{\epsilon}(\sigma)}$ to $\overline{\mathcal{M}(S)}$ by $\widehat{\overline{U_{\epsilon}(\sigma)}}$ which is a compact subset of $\overline{\mathcal{M}(S)}$. So the frame field $\mathbb{F}$ is in a compact subset of the extension of the frame bundle of the moduli space to the completion of the moduli space. 

Denote the unit disk bundle of the moduli space by $D\mathcal{M}(S)$. The short and relative length frame field provides a bundle extension for the tangent bundle of $\Teich(S)$ and $\mathcal{M}(S)$ to their WP  completions \cite{wolb}. Denote the extension of $D\mathcal{M}(S)$ to the completion by $D\overline{\mathcal{M}(S)}$. The vector $v$ projects to a vector in $D\widehat{\overline{U_{\epsilon}(\sigma)}}$.  Note that each of the functions $a_{\alpha}$, $b_{\alpha}$ and $c_{\beta}$, where $\alpha\in \sigma$ and $\beta \in \varsigma$ descends to a function on the compact set $K(\widehat{\overline{U_{\epsilon}(\sigma)}})\times D\widehat{\overline{U_{\epsilon}(\sigma)}}$. Therefore, each one of these functions is bounded on $\overline{U_{\epsilon}(\sigma)}$. Part (\ref{asympprod : ao(1)}) is proved.
\medskip

Since $v\in P_{\sigma}$, we have $v=\sum_{\alpha\in \sigma}a_{\alpha}\lambda_{\alpha}+b_{\alpha}J\lambda_{\alpha}$. Moreover since $w\in N_{\sigma}$, we have $w=\sum_{\beta\in \varsigma}c'_{\beta}\grad\ell_{\beta}$. Expanding $\langle v,w \rangle$ we get 
$$\langle v,w \rangle=\sum_{\alpha\in\sigma,\beta\in\varsigma}a_{\alpha}c'_{\beta}\langle \lambda_{\alpha},\grad\ell_{\beta} \rangle+b_{\alpha}c'_{\beta}\langle J\lambda_{\alpha},\grad\ell_{\beta} \rangle.$$
 By part (\ref{asympprod : ao(1)}) all of the coefficients in this sum are $O(1)$. Moreover by Proposition \ref{prop : wolpertmetric} all of the pairings in the sum above are $O(\epsilon)$. So each term in the sum above goes to $0$ as $\epsilon\to 0$. Part (\ref{asympprod : aprod}) is established.
\medskip

Now consider
 \begin{equation}\label{eq : vnormaldecomp}||v||^{2}=||v_{N}||^{2}+||v_{P}||^{2}+2\langle v_{N},v_{P}\rangle.\end{equation} 
 The vectors $\frac{v_{P}}{||v_{P}||}\in P_{\sigma}$ and $\frac{v_{N}}{||v_{N}||}\in N_{\sigma}$, and moreover $||\frac{v_{P}}{||v_{P}||}||^{2}=1$ and $||\frac{v_{N}}{||v_{N}||}||^{2}=1$. Then by part (\ref{asympprod : aprod}) for $\epsilon$ sufficiently small, 
 $$\langle \frac{v_{N}}{||v_{N}||}, \frac{v_{P}}{||v_{P}||} \rangle \leq \frac{1}{4}.$$
  If $||v_{N}||^{2}\geq ||v_{P}||^{2}$, then 
 $$|\langle v_{N},v_{P}\rangle|\leq \frac{1}{4}||v_{N}|| ||v_{P}||\leq \frac{1}{4}||v_{N}||^{2}.$$
 Substituing the above inequality into (\ref{eq : vnormaldecomp}) and dropping the non-negative term $||v_{P}||^{2}$ on the right-hand side, we get 
  $$||v||^{2}\geq ||v_{N}||^{2}-\frac{1}{2}||v_{N}||^{2}=\frac{1}{2}||v_{N}||^{2}.$$
   If $||v_{P}||^{2}\geq ||v_{N}||^{2}$, then by the exact same argument we get $||v||^{2}\geq \frac{1}{2}||v_{P}||^{2}$, which again implies that $||v||^{2}\geq \frac{1}{2}||v_{N}||^{2}$. This finishes the proof of part (\ref{asympprod : anorm}).
\end{proof}
 
\begin{lem}\label{lem : geodnstr} Given a multi-curve $\sigma$ and $T>0$, let $\zeta:[0,T]\to U_{\epsilon}(\sigma)$ be a WP geodesic segment parametrized by arc length. Let 
\begin{center}$\dot{\zeta}=\sum_{\alpha \in \sigma}a_{\alpha}\lambda_{\alpha}+b_{\alpha}J\lambda_{\alpha}+\sum_{\beta \in \varsigma} c_{\beta}\grad\ell_{\beta}$.\end{center}
 Then for every $\alpha \in \sigma$, we have that $a_{\alpha}\to 0$ and $b_{\alpha}\to 0$ as $\epsilon\to 0$. 
\end{lem}
\begin{proof}
We show that  for any $\alpha\in\sigma$ we have that $a_{\alpha}\to 0$ as $\epsilon\to 0$. The proof that $b_{\alpha}\to 0$ as $\epsilon\to 0$ is similar.

Define the constant 
\begin{eqnarray*}
\Lambda=\min\{||\lambda_{\alpha}(x)||^{2} : x\;\;\text{is in the $c_{0}-$thin part of the Teichm\"{u}ller space and}\; 
\alpha\in \mathcal{C}_{0}(S)\}.
\end{eqnarray*}
 Note that $\Lambda>0$ by the estimate for $||\lambda_{\alpha}||^{2}$ in Proposition \ref{prop : wolpertmetric}. To get a contradiction suppose that the above statement does not hold. Then there are a sequence $\epsilon_{n}\to 0$, a sequence of WP geodesic segments parametrized by arc-length $\zeta_{n}:[0,T]\to U_{\epsilon_{n}}(\sigma)$, and a sequence of times $t_{n}\in[0,T]$ so that for some $u>0$,
 $$a_{n}\geq \frac{3u}{\Lambda},$$
where $a_{n}$ is the coefficient of $\lambda_{\alpha}$ in the expansion of $\dot{\zeta}_{n}$ as in the statement of the lemma. After applying elements of the mapping class group (remarking), we may assume that the geodesic segments $\zeta_{n}$ are in a compact region of the completion of the Teichm\"{u}ller space. Then by the estimates in Proposition \ref{prop : wolpertmetric} for any $\alpha,\alpha'\in\sigma$ and $\beta\in\varsigma$ we have 
\begin{center}
$\langle \lambda_{\alpha}, \lambda_{\alpha'}\rangle=O(\epsilon_{n}^{3})$ ($\alpha\neq \alpha'$), $\langle \lambda_{\alpha}, J\lambda_{\alpha'}\rangle=0$ and $\langle \lambda_{\alpha},\grad\ell_{\beta} \rangle=O(\epsilon_{n}^{3/2})$. 
\end{center}
So we have that
  $$||\lambda_{\alpha}||^{2}a_{n}=\langle\lambda_{\alpha},\dot{\zeta}_{n}(t_{n}) \rangle+O(\epsilon_{n}^{3/2}).$$
   Thus, for all $n$ sufficiently large we have that
   \begin{equation}\label{eq : Fnterm}\langle \dot{\zeta}_{n}(t_{n}),\lambda_{\alpha}\rangle^{2}\geq 3u+O(\epsilon_{n}^{3/2})\geq 2u.\end{equation}
For each $n\in\mathbb{N}$, define the function 
$$F_{n}(t)=\langle \dot{\zeta}_{n}(t),\lambda_{\alpha}\rangle^{2}+\langle \dot{\zeta}_{n}(t),J\lambda_{\alpha}\rangle^{2}.$$
 By the formula at the end of $\S 5.2$ of \cite{uwpcurv}, $\frac{dF_{n}}{dt}=O(\ell_{\alpha}^{3/2})$, where the constant of the $O$ notation depends only on an upper bound for $\ell_{\alpha}$. For all $n$ sufficiently large we have $\zeta_{n}([0,T])\subset U_{\epsilon_{n}}(\sigma)$, so $\ell_{\alpha}(\zeta(t))\leq\epsilon_{n}$ for all $t\in[0,T]$. Moreover $\epsilon_{n}\to0$. Then for any $n$ large enough  
 $$\Big|\frac{dF_{n}}{dt}\Big|\leq \frac{u}{2T}.$$
Let $t\in[0,T]$.  Using the Fundamental theorem of calculus and the above inequality we get $|F_{n}(t)-F_{n}(t_{n})|\leq \frac{u}{2}$. By (\ref{eq : Fnterm}), $F_{n}(t_{n})\geq 2u$. Thus we may conclude that 
 \begin{equation}\label{eq : Fnlb}F_{n}(t)\geq \frac{3u}{2}\end{equation}
  for every $t\in[0,T]$. 

The formula from Proposition \ref{prop : cdlambda} for the WP covariant derivative and a straightforward calculation using properties of the Levi-Civita covariant derivative give us 
$$\ddot{\ell}_{\alpha}(\zeta_{n}(t))=2\langle \dot{\zeta}_{n}(t), \lambda_{\alpha}\rangle^{2} +6\langle\dot{\zeta}_{n}(t),J\lambda_{\alpha} \rangle^{2}+O(\epsilon_{n}^{3/2}).$$ 
Moreover, 
$$2\langle \dot{\zeta}_{n},\lambda_{\alpha}\rangle^{2}+6\langle \dot{\zeta}_{n},J\lambda_{\alpha}\rangle^{2} \geq 2F_{n}(t).$$
 Then, since $2F_{n}(t)\geq 3u$, by (\ref{eq : Fnlb}) for all $n$ sufficiently large, we obtain that
\begin{equation}\label{eq : l''}\ddot{\ell}_{\alpha}(\zeta_{n}(t))\geq u\end{equation} 
for every $t\in[0,T]$. 
\medskip

\begin{claim} For all $n$ sufficiently large we have that
\begin{equation}\label{eq : maxlzn}\max_{t\in [0,T]} \ell_{\alpha}(\zeta_{n}(t))\geq \frac{1}{16}uT^{2}.\end{equation}  
\end{claim}
Let $t_{\min,n}\in[0,T]$ be such that 
$$\ell_{\alpha}(\zeta_{n}(t_{\min,n}))=\min_{t\in [0,T]}\ell_{\alpha}(\zeta_{n}(t)).$$
 If there is a $t\in [0,T]$ such that 
$$\ell_{\alpha}(\zeta_{n}(t))-\ell_{\alpha}(\zeta_{n}(t_{\min,n}))\geq \frac{1}{16}uT^{2},$$
 then since $\ell_{\alpha}(\zeta_{n}(t_{\min,n}))\geq 0$, we get the desired lower bound in (\ref{eq : maxlzn}). Otherwise, 
 $$\ell_{\alpha}(\zeta_{n}(t))-\ell_{\alpha}(\zeta_{n}(t_{\min,n}))<\frac{1}{16}uT^{2}$$
  for every $t\in [0,T]$. Moreover, $\ell_{\alpha}(\zeta_{n}(t))\geq 0$ for every $t\in [0,T]$. Note that either the interval $[0,t_{\min,n}]$ or $[t_{\min,n},T]$ has length greater than or equal $\frac{T}{2}$. Then the Mean-value theorem for the interval with length at least $\frac{T}{2}$ implies that there is $t^{*}_{n}\in [0,T]$ such that $|\dot{\ell}_{\alpha}(t^{*}_{n})|<\frac{uT}{8}$. Let $t\in[0,T]$. Integrating both sides of the inequality (\ref{eq : l''}) from $t_{n}^{*}$ to $t$ and using the bound $|\dot{\ell}_{\alpha}(t^{*}_{n})|<\frac{uT}{8}$ we have
$$\ell_{\alpha}(\zeta_{n}(t))-\ell_{\alpha}(\zeta_{n}(t^{*}_{n}))\geq \frac{1}{2}u(t-t^{*}_{n})^{2}-\frac{uT}{8}(t-t^{*}_{n}).$$
Let $t'_{n}\in [0,T]$ be such that $|t'_{n}-t^{*}_{n}|=\frac{T}{2}$. At $t'_{n}$, the right-hand side of the above inequality is either $\frac{uT^{2}}{16}$ or $\frac{3uT^{2}}{16}$. Then since $\ell_{\alpha}(\zeta_{n}(t_{n}^{*}))\geq 0$, from the above inequality, we get 
$$\ell_{\alpha}(\zeta_{n}(t'_{n}))\geq \frac{uT^{2}}{16}.$$
 So we again have the lower bound (\ref{eq : maxlzn}). The proof of the claim is complete.
\medskip

The lower bound (\ref{eq : maxlzn}) for all $n$ sufficiently large contradicts the fact that $\zeta_{n}([0,T])$ (by the assumption that $\epsilon_{n}\to 0$) is a sequence of geodesic segments that converge to the $\sigma-$stratum where the $\alpha-$length-function is identically $0$. This contradiction completes the proof of that $a_{\alpha}\to 0$ as $\epsilon\to 0$.
\end{proof}

\begin{cor}\label{cor : geodP}
Let $\epsilon\leq c_{0}$. Let $\zeta$ be a WP geodesic segment in $U_{\epsilon}(\sigma)$. Then for $\epsilon$ sufficiently small, $\frac{1}{2}||\dot{\zeta}_{P}||^{2}\leq ||\dot{\zeta}||^{2}\leq ||\dot{\zeta}_{P}||^{2}$.
\end{cor}
\begin{proof}
 We have $||\dot{\zeta}||^{2}=||\dot{\zeta}_{P}||^{2}+||\dot{\zeta}_{N}||^{2}+2\langle \dot{\zeta}_{P}, \dot{\zeta}_{N} \rangle$. Let 
 \begin{center}$\dot{\zeta}=\sum_{\alpha \in \sigma}a_{\alpha}\lambda_{\alpha}+b_{\alpha}J\lambda_{\alpha}+\sum_{\beta \in \varsigma} c_{\beta}\grad\ell_{\beta},$\end{center} 
 then
\begin{center}$||\dot{\zeta}_{N}||^{2}=\sum_{\alpha,\alpha'\in\sigma}a_{\alpha}a_{\alpha'}\langle\lambda_{\alpha},\lambda_{\alpha'} \rangle+2a_{\alpha}b_{\alpha'}\langle\lambda_{\alpha},J\lambda_{\alpha'} \rangle+b_{\alpha}b_{\alpha'}\langle J\lambda_{\alpha},J\lambda_{\alpha'} \rangle .$\end{center}
 By Proposition \ref{prop : wolpertmetric} all of the inner products in the above sum are either $0$, $O(1)$ or $O(\epsilon^{3})$ where the constants of the $O$ notations depend only on $c_{0}$. Moreover by Lemma \ref{lem : geodnstr}, $a_{\alpha}$ and $b_{\alpha}$ go to $0$ as $\epsilon\to 0$. Thus $||\dot{\zeta}_{N}||^{2}\to 0$ as $\epsilon \to 0$. Furthermore, by Lemma \ref{lem : asympprod}(\ref{asympprod : aprod}), $\langle \dot{\zeta}_{P}, \dot{\zeta}_{N} \rangle\to 0$ as $\epsilon \to 0$. So for $\epsilon$ sufficiently small we have $\frac{1}{2}||\dot{\zeta}_{P}||^{2}\leq ||\dot{\zeta}||^{2}$. The inequality $||\dot{\zeta}||^{2}\leq ||\dot{\zeta}_{P}||^{2}$ is part (\ref{asympprod : anorm}) of Lemma \ref{lem : asympprod}.
\end{proof}

Let $\zeta$ be a smooth path and $V$ be a vector field along $\zeta$. We denote $\nabla_{\dot{\zeta}}V$ by $V'$ and $\nabla_{\dot{\zeta}}\nabla_{\dot{\zeta}}V$ by $V''$. The following lemma shows that the bundle $N_{\sigma}$ is almost parallel near the $\sigma-$stratum.

\begin{lem}\label{lem : N'}
Let $\epsilon\leq c_{0}$. Given a multi-curve $\sigma$ and $v>0$, let $\zeta:[0,T]\to U_{\epsilon}(\sigma)$ be a smooth curve and $V$ a vector field along $\zeta$ with $||V||^{2}\leq v$. Then 
$$\big| ||(V')_{N}||-||(V_{N})'||\big|\to 0$$
 as $\epsilon \to 0$. 
\end{lem}

\begin{proof}
Let 
\begin{eqnarray*}
V&=&\sum_{\alpha \in \sigma}a_{\alpha}\lambda_{\alpha}+b_{\alpha}J\lambda_{\alpha}+\sum_{\beta\in \varsigma}c_{\beta}\grad \ell_{\beta},\;\text{and}\\
 \dot{\zeta}&=&\sum_{\alpha \in \sigma}d_{\alpha}\lambda_{\alpha}+e_{\alpha}J\lambda_{\alpha}+\sum_{\beta \in \varsigma}f_{\beta}\grad\ell_{\beta}.
 \end{eqnarray*} 
 Expanding $V'_{N}$ and $(V_{N})'$, we get 
\begin{eqnarray*}
(V')_{N}-(V_{N})'&=&(\sum_{\beta\in \varsigma}c_{\beta}\nabla_{\dot{\zeta}}\grad\ell_{\beta}+\dot{c}_{\beta}\grad\ell_{\beta})_{N}=(\sum_{\beta\in \varsigma}c_{\beta}\nabla_{\dot{\zeta}}\grad\ell_{\beta})_{N}\\
&=&\sum_{\beta\in\varsigma,\alpha\in\sigma}\big(c_{\beta}d_{\alpha}\nabla_{\lambda_{\alpha}}\grad\ell_{\beta}+c_{\beta}e_{\alpha}\nabla_{J\lambda_{\alpha}}\grad\ell_{\beta}\big)_{N}\\
&+&\sum_{\beta,\beta'\in\varsigma}\big(c_{\beta}f_{\beta'}\nabla_{\grad\ell_{\beta'}}\grad\ell_{\beta}\big)_{N}
\end{eqnarray*}
By Proposition \ref{prop : wolpertcovderiv}, 
\begin{center}
$c_{\beta}d_{\alpha}\nabla_{\lambda_{\alpha}}\grad\ell_{\beta}=c_{\beta}d_{\alpha}O(\ell_{\alpha}^{1/2})$ and $c_{\beta}e_{\alpha}\nabla_{J\lambda_{\alpha}}\grad\ell_{\beta}=c_{\beta}e_{\alpha}O(\ell_{\alpha}^{1/2})$.\end{center}
 So each one of the terms in the first sum after the second equality above goes to $0$ as $\epsilon\to 0$. The vector fields $\grad\ell_{\beta}$ and $\grad\ell_{\beta'}$ are tangent to the stratum $\mathcal{S}(\sigma)$ so is $\nabla_{\grad\ell_{\beta'}}\grad\ell_{\beta}$. Therefore, by the continuity of the covariant derivatives as stated in Proposition \ref{prop : wolpertcovderiv}, 
$$(c_{\beta}f_{\beta'}\nabla_{\grad\ell_{\beta'}}\grad\ell_{\beta})_{N}\to 0$$
 as $\epsilon\to 0$. So we may conclude that all of the terms after the last equality above go to $0$ as $\epsilon \to 0$. So $ ||(V')_{N}-(V_{N})'||\to 0$ as $\epsilon \to 0$. Furthermore, 
$$\big| ||(V')_{N}||-||(V_{N})'||\big|\leq ||(V')_{N}-(V_{N})'||,$$
 hence $\big| ||(V')_{N}||-||(V_{N})'||\big|\to 0$ as $\epsilon\to 0$.
\end{proof}

\begin{proof}[Proof of Theorem \ref{thm : contnpt}] We start with some results in the general setting of a geodesically convex, negatively curved Riemannian manifold $\mathcal{X}$. Let $p\in \mathcal{X}$ and let $\pi(p)$ be the nearest point to $p$ on $g([0,T])$. Let $\zeta:[0,s_{1}]\to \mathcal{X}$ be the geodesic segment parametrized by arclength with $\zeta(0)=\pi(p)$ and $\zeta(s_{1})=p$. We have the following two situations:
\begin{enumerate}[(A)]
\item\label{spi(p) : I} $\pi(p)$ is in the interior of $g([0,T])$ or is an end point of $g([0,T])$ and is the nearest point to $p$ on a slightly longer geodesic segment containing $g([0,T])$ in its interior. 
\item\label{spi(p) : II} $\pi(p)$ is one of the end points $g(0)$ or $g(T)$ and is not the nearest point to $p$ on any geodesic segment containing $g([0,T])$ in its interior. 
\end{enumerate}
Suppose that (\ref{spi(p) : II}) holds. We claim that there is an open ball centered at $p$ which is mapped by the nearest point projection to one of the end points. Otherwise, there is a sequence of points $p_{n}\to p$ so that $\pi(p_{n})$ is the nearest point to $p_{n}$ on a slightly longer geodesic segment containing $g([0,T])$. But then continuity of the projection when points are projected to the interior of a geodesic segment (Proposition \ref{prop : proj map}(1)) implies that $\pi(p)$ is the nearest point to $p$ on a slightly longer geodesic segment containing $g([0,T])$. This contradicts the assumption (\ref{spi(p) : II}).

Now the existence of  a ball centered at $p$ which is mapped by the nearest point projection to one of the end points implies that $d\pi_{p}=0$. Thus $\frac{||d\pi_{p}(v)||^{2}}{||v||^{2}}<\delta$ for any $\delta\in[0,1)$.
\medskip

The rest of the proof is devoted to showing that $\frac{||d\pi_{p}(v)||^{2}}{||v||^{2}}<\delta$ for some $\delta\in[0,1)$ assuming (\ref{spi(p) : I}). First we reformulate the inequality in terms of Jacobi fields along the geodesic segments $\zeta$ connecting a point $p$ to its nearest point $\pi(p)$ on $g$. This reformulation will be convenient to work with. A vector field $J(s):=J(\zeta(s))$ along $\zeta(s)$ is a Jacobi field if it satisfies the Jacobi equation
\begin{equation}\label{eq : J}J'' +R(J,\dot{\zeta})\dot{\zeta} =0,\end{equation}
where $R(.,.).$ denotes the Riemann curvature operator. Moreover, $\dot{\zeta}$ denotes the derivative of $\zeta$ with respect to $s$. Also $J'=\nabla_{\dot{\zeta}}J$ and $J''=\nabla_{\dot{\zeta}}\nabla_{\dot{\zeta}}J$.  

Let us first characterize the map $d\pi : T\mathcal{X} \to Tg$ in terms of Jacobi fields. Given $v\in T_{p}\mathcal{X}$, let $\eta:[-\epsilon,\epsilon]\to \mathcal{X}$ be a smooth path passing through $p$ with $\eta(0)=p$ and $\dot{\eta}(0)=v$. Then the family of geodesics connecting each point $\eta(t)$ to  $\pi(\eta(t))$ the nearest point to $\eta(t)$ on $g$ defines a variation of geodesics 
$$u:\{(t,s): t\in[-\epsilon,+\epsilon], s\in[0,s_{1}]\}\to \mathcal{X},$$
where for each $t$, $u(t,s)$ is the geodesic connecting $\eta(t)$ to $\pi(\eta(t))$. Denote the geodesic connecting $p$ to $\pi(p)$ by $\zeta(s):=u(0,s)$, where $s\in[0,s_{1}]$. The vector field $\frac{\partial{u}}{\partial{t}}$ is a Jacobi field. Let $J(s):=\frac{\partial u}{\partial{t}}|_{(0,s)}$ for every $s\in [0,s_{1}]$, then $J(s_{1})=v$ and $J(0)=d\pi(v)$. 
\begin{lem} Let $J$ be a Jacobi field as above. Then we have
\begin{equation}\label{eq : J'}\frac{d}{ds}||J||^{2}|_{s=0}=0.\end{equation}
\end{lem}
\begin{proof}
A straightforward calculation shows that 
$$\frac{d}{ds}||J||^{2}|_{s=0}= 2\langle J'(0), J(0)\rangle.$$
 So we only need to verify that $\langle J(0),J'(0) \rangle=0$. Note that $\dot{\zeta}(s)=\frac{\partial{u}}{\partial{s}}|_{(0,s)}$ and $J(s)=\frac{\partial u}{\partial t}|_{(0,s)}$. Now we have 
\begin{eqnarray*}
\langle J(0),J'(0) \rangle&=&\langle \frac{\partial u}{\partial t}, \nabla_{\frac{\partial{u}}{\partial{s}}} \frac{\partial{u}}{\partial{t}} \rangle|_{s=0}=\langle \frac{\partial u}{\partial t}, \nabla_{\frac{\partial{u}}{\partial{t}}} \frac{\partial{u}}{\partial{s}} \rangle|_{s=0}\\
&=&-\langle \nabla_{\frac{\partial{u}}{\partial{t}}} \frac{\partial{u}}{\partial{t}},\frac{\partial {u}}{\partial {s}} \rangle|_{s=0}
=-\langle \nabla _{f\dot{g}} f\dot{g},\frac{\partial{u}}{\partial{s}} \rangle|_{s=0}\\
&=&- \langle f\dot{g}(f)\dot{g},\frac{\partial{u}}{\partial{s}} \rangle|_{s=0}=0.
\end{eqnarray*}
The second equality above follows because $\nabla$ is torsion free and $[\frac{\partial u}{\partial s},\frac{\partial u}{\partial t}]=0$. To get the third equality note that, by the compatibility of the Levi-Civita covariant derivative and the Riemannian metric we have: 
$$\frac{d}{dt}\langle\frac{\partial u}{\partial t},\frac{\partial u}{\partial s} \rangle=\langle \nabla_{\frac{\partial{u}}{\partial{t}}} \frac{\partial{u}}{\partial{t}},\frac{\partial {u}}{\partial {s}} \rangle+\langle \nabla_{\frac{\partial{u}}{\partial{t}}} \frac{\partial{u}}{\partial{t}},\frac{\partial {u}}{\partial {s}} \rangle.$$
 We also have that $\langle \frac{\partial{u}}{\partial{t}},\frac{\partial{u}}{\partial{s}} \rangle|_{s=0}=0$ for every $t$, and hence $\frac{d}{dt}\langle\frac{\partial u}{\partial t},\frac{\partial u}{\partial s} \rangle=0$. Then the equality follows from the above equality. In the fourth equality we replace $\frac{\partial{u}}{{\partial{t}}}|_{s=0}$ by $f\dot{g}$, where $f:[-\epsilon,+ \epsilon]\to \mathbb{R}$ is a function. The fifth equality follows from a straightforward calculation using the fact that $\nabla_{\dot{g}}\dot{g}=0$ for the geodesic $g$. The last equality holds because $\frac{\partial u}{\partial t}|_{s=0}=f \dot{g}$ and $\langle \frac{\partial{u}}{\partial{t}},\frac{\partial{u}}{\partial{s}} \rangle|_{s=0}=0$ for all $t$. 
\end{proof}

Using the Jacobi equation (\ref{eq : J}) for any $s\in[0,s_{1}]$ we have
\begin{equation}\label{eq : J''}\frac{1}{2}\frac{d^{2}}{ds^{2}}||J||^{2}=||J'||^{2}-\kappa(s)|\dot{\zeta} \wedge J|^{2},\end{equation}
where $\kappa(s)=\kappa(\dot{\zeta}(s),J(s))$ is the sectional curvature of the span of $\dot{\zeta}(s)$ and $J(s)$ and 
$$|\dot{\zeta} \wedge J|^{2}=||J||^{2}||\dot{\zeta}||^{2}-\langle J,\dot{\zeta} \rangle^{2}.$$ 

\begin{lem}\label{lem : obsJ} Given a Jacobi field $J$ we have
\begin{enumerate}[(i)]
\item $\frac{d}{ds}||J(s)||^{2}\geq 0$ for every $s\in [0,s_{1}]$.
\item $||J(s)||^{2}$ is non-decreasing on the interval $[0,s_{1}]$.
\item $\langle J(s),\dot{\zeta}(s)\rangle$ is a linear function. 
\end{enumerate} 
\end{lem}
\begin{proof}
 By  (\ref{eq : J''}) and the fact that $\kappa\leq 0$ we have that 
 $$\frac{d^{2}}{ds^{2}}||J||^{2}\geq 0.$$
  Thus $\frac{d}{ds}||J||^{2}$ is a non-decreasing function of $s$. Furthermore, by (\ref{eq : J'}), 
  $$\frac{d}{ds}||J||^{2}|_{s=0}=0.$$
   Therefore we have that $\frac{d}{ds}||J||^{2}\geq 0$. Part (i) is proved. Part (ii) follows from part (i). We proceed to prove part (iii). A straightforward calculation using the Jacobi equation (\ref{eq : J}) gives us 
$$\frac{d^{2}}{ds^{2}}\langle J,\dot{\zeta}\rangle=\langle R(J,\dot{\zeta})\dot{\zeta},\dot{\zeta} \rangle=0.$$
This implies that $\langle J(s),\dot{\zeta}(s)\rangle$ is linear, because its second derivative is identically $0$.  Part (iii) is proved.
\end{proof}
Given $b\geq 0$, denote the closed $b-$neighborhood of $g([0,T])\subset \mathcal{X}$ by $\overline{\mathcal{N}}_{b}(g([0,T]))$. Denote the boundary of $\overline{\mathcal{N}}_{b}(g([0,T]))$ by $\partial{\overline{\mathcal{N}}}_{b}\mathcal(g([0,T]))$. Since $\mathcal{X}$ is negatively curved, the nearest point projection map $\pi:\mathcal{X}\to \partial{\overline{\mathcal{N}}}_{b}\mathcal(g([0,T]))$ is 1-Lipschitz. This follows from the growth of the norm of Jacobi fields given in the proof of part (\ref{pmetric : distinc}) of Theorem \ref{thm : pullbackm}. Furthermore, the projection of a point $p\in \mathcal{X}$ to $g([0,T])$ may be obtained by first projecting $p$ on $\partial{\overline{\mathcal{N}}}_{b}(g([0,T]))$ and then projecting the projected point to $g([0,T])$. So we only need to prove the strict contraction of the nearest point projection map for the points of $\partial{\overline{\mathcal{N}}}_{b}(g([0,T]))$.
 
Now we return to the setting of the Teichm\"{u}ller space equipped with the WP metric. Recall that $g([0,T])\subset U_{\epsilon,\bar{\epsilon}}(\sigma)$. Let $p\in \partial{\overline{\mathcal{N}}}_{b}(g([0,T]))$ and let $\zeta$ be the geodesic segment connecting $p$ to $\pi(p)$. Given a unit vector $v\in T_{p}\mathcal{X}$ ($||v||^{2}=1$) as we saw earlier there is a Jacobi field $J$ along $\zeta$ with $J(b)=v$ and $J(0)=d\pi(v)$. Then $\frac{||d\pi_{p}(v)||^{2}}{||v||^{2}}= ||J(0)||^{2}$. In Theorem \ref{thm : normJbd} below we show that there is a $\delta>0$, so that $ ||J(0)||^{2}\leq\delta<1$. This completes the proof of the theorem.
\end{proof}

\begin{thm}\label{thm : normJbd}
Given $\bar{\epsilon}>0$ and $b,b'>0$ sufficiently small with $b>b'$, there is a constant $\delta=\delta(\bar{\epsilon},b,b')$ with the following property. Let $\epsilon>0$ be sufficiently small. Let $\sigma$ be a multi-curve such that $S\backslash\sigma$ is a large subsurface. Let $g([0,T])\subset U_{\epsilon,\bar{\epsilon}}(\sigma)$ be a geodesic segment. Let $p\in\partial{\overline{\mathcal{N}}}_{b}(g([0,T]))$ and suppose that the nearest point $\pi(p)$ to $p$ on $g([0,T])$ be as in (\ref{spi(p) : I}) in the proof  of Theorem \ref{thm : contnpt}. Let $v\in T_{p}\mathcal{X}$ and let $\zeta$ be a geodesic connecting $p$ to $\pi(p)$. Finally, let $J$ be a Jacobi field we setup in the proof of Theorem \ref{thm : contnpt}. Then 
 \begin{equation}\label{eq : Js<Jb}||J(b')||^{2}\leq\delta.\end{equation}
Moreover $\delta$ decreases as $b'\to 0$. 
\end{thm}

\begin{proof}
 First we sketch the proof. Suppose that a neighborhood of a stratum $\mathcal{S}(\hat{\sigma})$, where $\hat{\sigma}\subseteq\sigma$, is foliated with totally geodesic leaves with negative sectional curvatures bounded away from $0$. If both the geodesic $\zeta$ connecting $p$ to $\pi(p)$ and $J$ are tangent to the leaves of the foliation over an interval $E$, then $\kappa(J,\dot{\zeta})$ is negative over $E$ and we obtain a negative upper bound for the second term on the right-hand side of (\ref{eq : J''}). Otherwise, $J$ would vary over $E$, so that we may obtain a lower bound for $||J'||$, and therefore the first term on the right-hand side of (\ref{eq : J''}). Having these bounds, integrating both sides of the equation (\ref{eq : J''}) we obtain an upper bound for $||J(b')||^{2}$.
\medskip

We do not quite have the above picture in the WP setting. However a modification of the above argument, as is outlined below, would give us the bound. As we will see in Claims \ref{claim : zinethicklbd} and \ref{claim : intervalE} below, there is an $\epsilon'>0$, so that for any $u<\epsilon'$ there exist an interval $E$ and a multi-curve $\hat{\sigma}\subset\sigma$ of definite length such that $\zeta(t)\in U_{u,\epsilon'}(\hat{\sigma})$. Moreover, instead of the totally geodesic foliation with negatively curved leaves near the stratum  we have the decomposition of the tangent bundle of the Teichm\"{u}ller space 
$$T\Teich(S)=P_{\hat{\sigma}}+N_{\hat{\sigma}}.$$
 If the $P_{\hat{\sigma}}$ component of $J(s)$ has a definite length over $E$, then using compactness arguments we establish a negative upper bound for the sectional curvature of the span of $\dot{\zeta}(s)$ and $J(s)$ over a subinterval of $E$. This provides a negative upper bound for the second term on the right-hand side of (\ref{eq : J''}). Note that the sectional curvatures are bounded away from $0$ in the thick part of the $\hat{\sigma}-$stratum ($\hat{\sigma}$ is a non-separating multi-curve) and $P_{\hat{\sigma}}$ is almost tangent to the level manifolds of the functions $(\ell_{\alpha}^{1/2})_{\alpha \in \hat{\sigma}}$ which define a foliation in a neighborhood of $\mathcal{S}(\hat{\sigma})$. This is Case 1 below.
 
  If the $P_{\hat{\sigma}}$ component of $J(s)$ does not have a definite length over $E$, then since $J(0)\in P_{\hat{\sigma},\zeta(0)}$ (the normal component of $J(0)$ is $0$), the normal component of $J$ varies over $E$ which gives rise to a lower bound for the integral of $||(J_{N})'||$ on a subinterval of $E$. Furthermore, by Lemma \ref{lem : N'} the bundle $N_{\hat{\sigma}}$ is almost parallel, so we may obtain a lower bound for the integral of  $||J'||^{2}$. This provides a lower bound for $||J'||^{2}$ and consequently the first term on the right-hand side of  (\ref{eq : J''}). Having these bounds integrating both sides of the equation (\ref{eq : J''}) we get the desired bound for $||J(b')||^{2}$. This is Case 2 below.
\medskip

We proceed with the proof of the theorem following the above sketch. The following two claims provide us a thick-thin decomposition for a portion of the geodesic segment $\zeta$.

\begin{claim}\label{claim : zinethicklbd}
For any $\bar{\epsilon}>0$, there are $h,\epsilon'>0$ with the following properties. Suppose that $\zeta$ is a geodesic segment in the $h$ neighborhood of the $\bar{\epsilon}-$thick part of $\mathcal{S}(\sigma)$. Then $\ell_{\beta}(\zeta(t))>\epsilon'$ for any $\beta\notin\sigma$.
\end{claim}
Passing to the quotient and using the fact that there are finitely many disjoint strata in the completion of the moduli space $\overline{\mathcal{M}(S)}$, there is a lower bound for the distance between $\overline{\mathcal{S}(\sigma)}$ and $\overline{\mathcal{S}(\tau)}$ for any two multi-curves $\sigma$ and $\tau$ with $\sigma\pitchfork\tau$. Moreover the $\bar{\epsilon}-$thick part of the $\sigma-$stratum is compact and does not contain any point in $\overline{\mathcal{S}(\sigma)}\backslash\mathcal{S}(\sigma)$. So there is a lower bound for the distance between the $\bar{\epsilon}-$thick part of the $\sigma-$stratum and any point in a $\tau-$stratum, where $\sigma\subsetneq\tau$. Therefore, there is a lower bound for the distance between the $\bar{\epsilon}-$thick part of $\mathcal{S}(\sigma)$ and any $\tau-$stratum with $\tau\pitchfork\sigma$ or $\sigma\subsetneq\tau$ only depending on $\bar{\epsilon}$. This means that there is an $h>0$ such that the $h-$neighborhood of the $\bar{\epsilon}-$thick part of $\mathcal{S}(\sigma)$ only intersects the strata of multi-curves $\sigma'\subseteq \sigma$. Now suppose that $\zeta$ is in the $h$ neighborhood of the $\bar{\epsilon}-$thick part of $\mathcal{S}(\sigma)$. Thena compactness argument shows that there is a lower bound $\epsilon'>\epsilon$ for the length of every curve $\beta\notin\sigma$ along $\zeta$. The claim is proved.
\medskip

In the rest of the proof we assume that $b\leq\frac{h}{2}$. 

\begin{claim}\label{claim : intervalE}
There is an $e>0$ with the following property. Let $\zeta:[0,b]\to\Teich(S)$ be the geodesic segment connecting a point $p\in\partial{\overline{\mathcal{N}}_{b}(g([0,T]))}$ to the nearest point to $p$ on $g([0,T])$. For any $u<\epsilon'$, there is a multi-curve $\hat{\sigma}\subseteq\sigma$ (possibly empty) and an interval $E\subseteq [0,b]$ with $|E|= e$ such that $\zeta(t)\in U_{u}(\hat{\sigma})$ for any $t\in E$. 
\end{claim}

Since $g([0,T])\subset U_{\epsilon,\bar{\epsilon}}(\sigma)$, the geodesic segment $g([0,T])$ converges into the $\bar{\epsilon}-$thick part of $\mathcal{S}(\sigma)$ as $\epsilon\to 0$. Thus for $\epsilon$ sufficiently small $g([0,T])$ is contained in the $\frac{h}{2}$ neighborhood of the $\bar{\epsilon}-$thick part of $\mathcal{S}(\sigma)$. $\overline{\mathcal{N}}_{h}(g([0,T]))$ is geodesically convex, and the points $\zeta(0),\zeta(b)$ lie in $\overline{\mathcal{N}}_{h}(g([0,T]))$. So $\zeta([0,b])\subset\overline{\mathcal{N}}_{h}(g([0,T]))$ and therefore $\zeta([0,b])$ is in the $h$ neighborhood of the $\bar{\epsilon}-$thick part of $\mathcal{S}(\sigma)$. 

Given $t\in [0,b]$, let $\sigma_{t}$ be the maximal subset of $\sigma$ (possibly empty) such that $\ell_{\alpha}(\zeta(t))< u$ for any $\alpha\in\sigma_{t}$. Since $\zeta([0,b])$ is in the $h$ neighborhood of the $\bar{\epsilon}-$thick part of $\mathcal{S}(\sigma)$, by Claim \ref{claim : zinethicklbd}, we have that $\ell_{\alpha}(\zeta(t))>\epsilon'$ for any $\alpha\notin\sigma$. Then since $u\leq\epsilon'$, we may conclude that $\sigma_{t}\subseteq\sigma$.
  
By the convexity of length-functions along WP geodesics (Theorem \ref{thm : convlf}), for any $\sigma'\subseteq \sigma$ the set of times $t\in[0,b]$ such that $\sigma_{t}=\sigma'$ is an interval. Moreover, the number of multi-curves $\sigma'\subseteq \sigma$ is $2^{|\sigma|}\leq 2^{\xi(S)}$ where $\xi(S)=3g-3+n$. By these two observations there is a possibly empty multi-curve $\hat{\sigma}\subseteq\sigma$ and an interval $E$ of length at least $e =\frac{b}{2^{\xi(S)}}$ such that $\sigma_{t}=\hat{\sigma}$ for every $t\in E$. Note that $e$ only depends on $b$ and $\bar{\epsilon}$. The proof of the claim is complete.
 \medskip

Any vector in $T_{p}\mathcal{X}$ has a decomposition into a component tangent to $\zeta$ and a component orthogonal to $\zeta$. We are in the situation that $\langle \dot{\zeta},\dot{g}\rangle|_{\pi(p)}=0$ and $d\pi_{p}:T_{p}\mathcal{X}\to T_{\pi(p)}g$ is a linear map; see (\ref{spi(p) : I}). The linear map $d\pi$ maps any vector tangent to $\zeta$ to $0\in T_{\pi(p)}g$. Thus by the linearity of $d\pi$ we only need to prove (\ref{eq : Js<Jb}) for vectors $v$ with $\langle v,\dot{\zeta}(b)\rangle=0$. Let $v$ be such a vector, and suppose that $||v||^{2}=1$. Let $J$ be a Jacobi field with $J(0)=d\pi(v)$ and $J(b)=v$. Then since $\langle\dot{\zeta}(b),J(b) \rangle=0$ and $\langle \dot{\zeta}(0),J(0)\rangle=0$, by Lemma \ref{lem : obsJ}(ii), for every $s\in[0,b]$ we have 
\begin{equation}\label{eq : Jz'orth}
\langle J(s), \dot{\zeta}(s)\rangle= 0.
\end{equation}
Let $s\in[0,b]$. By Lemma \ref{lem : obsJ}(i), $||J(s)||^{2}$ is a non-decreasing function of $s$. Then since $||J(b)||^{2}=1$, we have $||J(s)||^{2}\leq 1$. Fix a constant $\omega_{0}<1$. If $||J(0)||^{2}\leq \omega_{0}$, then for $\delta=\frac{1+\omega_0}{2}$ and any $b,b'$ sufficiently small the bound (\ref{eq : Js<Jb}) holds by continuity of $||J(s)||^{2}$. Moreover, the fact that $||J(s)||^{2}$ is non-decreasing in $s$ implies that $\delta$ is non-increasing as $b'\to 0$. So in the rest of the proof we will assume that $||J(0)||^{2}\geq \omega_{0}$. Then using the fact that $||J(s)||^{2}$ is non-decreasing we have $||J(s)||^{2}\geq \omega_{0}$. We record the upper and lower bounds for $||J(s)||^{2}$,
\begin{equation}\label{eq : normJbd}
\omega_{0}\leq ||J(s)||^{2}\leq 1. 
\end{equation}

The geodesic $g([0,T])$ is contained in $U_{\epsilon}(\sigma)$, so by Corollary \ref{cor : geodP}, $||\dot{g}_{P}||^{2}\geq \frac{1}{2}||\dot{g}||^{2}$. Then since $J(0)=||J(0)||\dot{g}$, by (\ref{eq : normJbd}), we have $||J_{P}(0)||^{2}\geq \frac{\omega_{0}}{2}$. 

For $u$ sufficiently small, let $E$ be the interval from Claim \ref{claim : intervalE}. We will consider the following two cases depending on the behavior of the function $||J_{P}(s)||^{2}$ on $E$, and in each case verify that the inequality (\ref{eq : Js<Jb}) holds.
 \medskip

\noindent\textbf{Case 1:} We have
 \begin{equation}\label{eq : c1}||J_{P}(s)||^{2}\geq \frac{\omega_{0}}{4}\end{equation}
  for every $s\in E$.
  \medskip
  
  We prove that the inequality (\ref{eq : Js<Jb}) holds using the following lemma:
   \begin{lem} \label{lem : seccurvpub}
 There exist a closed subinterval $E'\subseteq E$ and  $k_{0}<0$ such that for every $s\in E'$ we have $\kappa(\dot{\zeta}(s),J(s))\leq k_{0}$.
\end{lem}
 Let the interval $E'=[c,d]$ and $k_{0}<0$ be from Lemma \ref{lem : seccurvpub}. Let $s\in [c,d]$. Then by the lemma we have that
$$\kappa(s)=\kappa(\dot{\zeta}(s),J(s))\leq k_{0}.$$
 Moreover, $[c,d]\subset [0,b]$, so by (\ref{eq : Jz'orth}), $\langle J(s), \dot{\zeta}(s)\rangle = 0$, and by (\ref{eq : normJbd}), $||J(s)||^{2}\geq\omega_{0}$. Thus we have that
 $$|J \wedge \dot{\zeta}|^{2} =||J||^{2}||\dot{\zeta}||^{2}-\langle J,\dot{\zeta}\rangle\geq \omega_{0}.$$
  Using the above two inequalities we obtain
 \begin{equation}\label{eq : -k|J^z.|lb}-\kappa(s)|J \wedge \dot{\zeta}|^{2}\geq -k_{0}\omega_{0}>0.\end{equation}
Moreover by (\ref{eq : J''}),
 $$\frac{1}{2}\frac{d^{2}}{ds^{2}}||J(s)||^{2}= ||J'(s)||^{2}-\kappa(s)|J(s)\wedge \dot{\zeta}(s)|^{2},$$
  then by the inequality (\ref{eq : -k|J^z.|lb}) we have
$$\frac{d^{2}}{ds^{2}}||J(s)||^{2}\geq -2k_{0}\omega_{0}>0.$$ 
Suppose that $b'\leq |E'|$, then $b'\leq d$. Integrating both sides of the above inequality over the interval $[b',s]$ ($s>b'$) twice, we have 
$$\frac{d}{ds}||J(s)||^{2}-\frac{d}{ds}||J(b')||^{2}\geq -2k_{0}\omega_{0}(s-b').$$ 
Moreover, by Lemma \ref{lem : obsJ}(iii), $\frac{d}{ds}||J(s)||^{2}\geq 0$. So we get 
$$\frac{d}{ds}||J(s)||^{2}\geq -2k_{0}\omega_{0}(s-b').$$
 Now integrating both sides of the above inequality from $b'$ to $d$ we get   
\begin{eqnarray}\label{eq : normJ(d)}
||J(d)||^{2}&\geq& ||J(b')||^{2}-k_{0}\omega_{0}(s-b')^{2}\big|_{b'}^{d}\\
&=& ||J(b')||^{2}-k_{0}\omega_{0}(d-b')^{2}.\nonumber
\end{eqnarray}
Moreover, by Lemma \ref{lem : obsJ}(i), $||J(s)||^{2}$ is non-decreasing, so 
\begin{center}$||J(b)||^{2}\geq ||J(d)||^{2}$ and $||J(b')||^{2}\geq ||J(0)||^{2}$.\end{center} 
Then from (\ref{eq : normJ(d)}) we may deduce that
$$1=||J(b)||^{2}\geq ||J(0)||^{2}-2k_{0}\omega_{0}(b-b')^{2}.$$
 Consequently, for $\delta=1+k_{0}\omega_{0}(b-b')^{2}$, (\ref{eq : Js<Jb}) holds. Here by choosing $b$ sufficiently small we can guarantee that $\delta>0$. Furthermore, note that since $k_{0}<0$, the number $\delta$ is less than $1$ and decreases as $b'\to 0$.
\medskip

We finish Case 1 by proving Lemma~\ref{lem : seccurvpub}. 

\begin{proof}[Proof of Lemma \ref{lem : seccurvpub}]
 The proof of the lemma is by contradiction. Suppose that the lemma does not hold. Then there exist a sequence of points $p_{n}\in\partial{\overline{\mathcal{N}}_{b}(g([0,T]))}$ and geodesic segments $\zeta_{n}:[0,b]\to\Teich(S)$ connecting $p_{n}$ to the nearest point to $p_{n}$ on $g([0,T])$, and intervals $E_{n}\subseteq[0,b]$ as in Claim \ref{claim : intervalE}, so that for any sequence of times $s_{n}\in I_{n}$ we have
 \begin{equation}\label{eq : curv zeta J at sn}\kappa(\dot{\zeta}_{n}(s_{n}),J_{n}(s_{n}))\to 0\end{equation}
  as $n\to \infty$, where $J_{n}$ is a Jacobi field along the geodesic segment $\zeta_{n}$ which satisfies (\ref{eq : Jz'orth}), (\ref{eq : normJbd}) and (\ref{eq : c1}). 
  
Note that $|E_{n}|\geq\frac{b}{2^{\xi(S)}}$ for each $n\in\mathbb{N}$. So after trimming the intervals $E_{n}$ and reparametrizing the geodesic segments $\zeta_{n}|_{E_{n}}$ we may assume that $E_{n}\equiv E$.
  
  \begin{claim}\label{claim : intervalE'}
There exist $u'>0$, a sequence $u_{n}\to 0$ as $n\to\infty$, a possibly empty multi-curve $\hat{\tau}\subseteq\hat{\sigma}$ and a closed subinterval $E'\subseteq E$ so that $\zeta_{n}(t)\in U_{u_{n},u'}(\hat{\tau})$ for any $t\in E'$. 
\end{claim}

The proof is by contradiction.  Apply Theorem \ref{thm : geodlimit} to the geodesic segments $\zeta_{n}$ and let the partition $t_{0}<....<t_{k+1}$, multi-curves $\sigma_{i}$ for $i=0,...,k+1$, $\hat{\tau}$ and the piece-wise geodesic $\hat{\zeta}$ be as in the theorem. Let $0\leq j\leq k$ and let $E'$ be a closed subinterval of $(t_{j},t_{j+1})$. Then 
  \begin{itemize}
  \item there is $u''>0$, so that for any $t\in E'$ and any $\alpha\in\hat{\sigma}-\hat{\tau}$ we have $\ell_{\alpha}(\hat{\zeta}(t))>2u''$,
  \item  for any $\alpha\in\hat{\tau}$ and any $t\in E'$, $\ell_{\alpha}(\hat{\zeta}(t))=0$.
  \end{itemize}
  Let the elements $\psi_{n}$ of mapping class group be as in Theorem \ref{thm : geodlimit}. The geodesic segments $\zeta_{n}$ are in $\overline{\mathcal{N}}_{b}(g([0,T]))$, which is a compact subset of $\overline{\Teich(S)}$. These two facts imply that $\psi_{n}\equiv I$. Moreover by Claim \ref{claim : intervalE}, $\ell_{\alpha}(\zeta(t))>\epsilon'$ for any $\alpha\notin\hat{\sigma}$, so for each $i=1,...,k$ we have $\sigma_{i}\subseteq\hat{\sigma}$. Similarly $\hat{\tau}\subseteq\hat{\sigma}$. For $i=1,...,k$ and any $n\in\mathbb{N}$, let $\mathcal{T}_{i,n}$ and $\varphi_{i,n}$ be as in Theorem \ref{thm : geodlimit}. Recall that each $\mathcal{T}_{i,n}\in\tw(\hat{\sigma})$. Thus each $\varphi_{i,n}$ is the composition of $I$ and powers of positive Dehn twists about curves in $\hat{\sigma}$. So applying each $\varphi_{i,n}$ to a point in the Teichm\"{u}ller space does not change the isotopy class and length of any curve $\alpha\in\hat{\sigma}$. Let $E'\subset (t_{j},t_{j+1})$ be as above, by Theorem \ref{thm : geodlimit}(3), for any $t\in E' $ we have $\varphi_{j,n}(\zeta_{n}(t))\to\hat{\zeta}(t)$ as $n\to\infty$. Then since $\varphi_{j,n}$ does not change the length and isotopy class of curves in $\hat{\tau}$, it follows from the continuity of length-functions and the two bullets above, respectively, that 
\begin{itemize}
\item$\ell_{\alpha}(\zeta_{n}(t))> u''$ for any $\alpha\in\hat{\tau}$ and $t\in E'$,
\item there is a sequence $u_{n}\to 0$ as $n\to\infty$, so that $\ell_{\alpha}(\zeta_{n}(t))<u_{n}$ for any $\alpha\in\hat{\tau}$ and $t\in E'$. 
\end{itemize}
Also since $E'\subseteq E$, as we mentioned before, by Claim \ref{claim : intervalE}, we have $\ell_{\alpha}(\zeta_{n}(t))> \epsilon'$, for any $\alpha\notin\hat{\sigma}$, and $t\in E'$. 

Thus the claim holds for the interval $E'\subseteq E$, the multi-curve $\hat{\tau}$, the constant $u'=\min\{\epsilon',u''\}$, and the sequence $u_{n}$. 
 \medskip

The geodesic segments $\zeta_{n}(E')$ are in $\overline{\mathcal{N}}_{b}(g([0,T]))$, which is a compact subset of $\overline{\Teich(S)}$. Then by Claim \ref{claim : intervalE'} after possibly passing to a subsequence the geodesic segments $\zeta_{n}(E')$ converge to the $u'-$thick part of $\mathcal{S}(\hat{\tau})$.  

First suppose that $\hat{\tau}\neq\emptyset$.  Let $s_{n}\in E'$ and the points $\zeta_{n}(s_{n})$ converge to a point $p$ as $n\to\infty$. The vectors $\dot{\zeta}_{n}(s_{n})$ converge to a vector $v\in T_{p}\Teich(S)$ with $||v||^{2}=1$ and the vectors $J_{n}(s_{n})$ converge to a vector $w\in T_{p}\Teich(S)$ with $||w||^{2}\geq \omega_{0}$. Moreover, $\langle v,w\rangle=0$. So $\spn(v,w)$ is a nondegenerate plane at $p$ in the $u'-$thick part of the Teichm\"{u}ller space. There is an upper bound $k_{1}<0$ for the sectional curvatures in this region. Thus for $n$ sufficiently large, 
$$\kappa(\dot{\zeta}_{n}(s_{n}),J_{n}(s_{n}))<\frac{k_{1}}{2}<0.$$
But this contradicts the assumption that (\ref{eq : curv zeta J at sn}) holds.
\begin{remark}
The upper bound for the sectional curvatures in the thick part of the moduli space follows from the fact that the Weil-Petersson metric is defined on the moduli space and the thick part of the moduli space is compact. 
\end{remark}

 In the rest of the proof of the lemma we assume that $\hat{\tau}\neq\emptyset$. Since $\hat{\tau}$ is a non-separating multi-curve, the stratum $\mathcal{S}(\hat{\tau})$ is not a product of lower complexity strata. Then all of the sectional curvatures of the stratum $\mathcal{S}(\hat{\tau})$ are bounded above by some $k_{2}<0$ in the $u'-$thick part of the stratum. 
 
 Let $s_{n}\in E'$ and the points $\zeta_{n}(s_{n})$ converge to the point $p$ as $n\to\infty$. In what follows for simplicity of notation we often drop reference to the point $\zeta_{n}(s_{n})$. For example, we write $\lambda_{\alpha}$ instead of $\lambda_{\alpha}(\zeta_{n}(s_{n}))$. Let
  \begin{eqnarray*}
 \dot{\zeta}_{n}(s_{n})&=&\sum_{\alpha \in \hat{\sigma}}a_{\alpha,n}\lambda_{\alpha}+b_{\alpha,n}J\lambda_{\alpha}+\sum_{\beta \in \varsigma}c_{\beta,n}\grad\ell_{\beta},\;\text{and}\\
  J_{n}(s_{n})&=&\sum_{\alpha \in \hat{\sigma}}d_{\alpha,n}\lambda_{\alpha}+e_{\alpha,n}J\lambda_{\alpha}+\sum_{\beta\in \varsigma}f_{\beta,n}\grad \ell_{\beta}. 
 \end{eqnarray*}
 Since $\zeta_{n}(s_{n})\in U_{u_{n}}(\hat{\tau})$ and $||\dot{\zeta}_{n}||^{2}=1$, by Lemma \ref{lem : asympprod}(\ref{asympprod : ao(1)}), the coefficients $a_{\alpha,n}$, $b_{\alpha,n}$ and $c_{\beta,n}$ are $O(1)$. Similarly, since $||J_{n}||^{2}\leq 1$, the coefficients $d_{\alpha,n},e_{\alpha,n}$ and $f_{\beta,n}$ are $O(1)$, where the constants of $O$ notations depend only on an upper bound for $u'$ and $u_{n}$'s.
 
We have the following bounds for the terms in the expansion of $\langle R(\dot{\zeta}_{n},J_{n})J_{n},\dot{\zeta}_{n} \rangle$: 
\begin{itemize}
\item For any $\alpha\in\hat{\tau}$, the four terms 
\begin{itemize}
\item $a_{\alpha,n}^{2}e_{\alpha,n}^{2}\langle R(\lambda_{\alpha},J\lambda_{\alpha})J\lambda_{\alpha},\lambda_{\alpha}\rangle$, 
\item $b_{\alpha,n}^{2}d_{\alpha,n}^{2}\langle R(J\lambda_{\alpha},\lambda_{\alpha})\lambda_{\alpha},J\lambda_{\alpha}\rangle$,
\item $a_{\alpha,n}e_{\alpha,n}d_{\alpha,n}b_{\alpha,n}\langle R(\lambda_{\alpha},J\lambda_{\alpha})\lambda_{\alpha},J\lambda_{\alpha}\rangle$, 
\item $a_{\alpha,n}e_{\alpha,n}d_{\alpha,n}b_{\alpha,n}\langle R(J\lambda_{\alpha},\lambda_{\alpha})J\lambda_{\alpha},\lambda_{\alpha}\rangle$, 
 \end{itemize} 
 add up to a non-positive number for all $n$ sufficiently large.
\item For any $\alpha\in\hat{\tau}$, any term which is a multiple of $\langle R((J)\lambda_{\alpha},.).,.\rangle$ except the ones in the first bullet above is either $0$ or has arbitrarily small absolute value for all $n$ sufficiently large.
\end{itemize}
{\bf Proof of the first bullet point:} For each $n$, $\langle R(\lambda_{\alpha},J\lambda_{\alpha})J\lambda_{\alpha},\lambda_{\alpha}\rangle$ is evaluated at the point $\zeta_{n}(s_{n})$ which converges to the point $p$ in the $u'-$thick part of $\mathcal{S}(\hat{\tau})$. Then by the limit of the diagonal terms in Proposition \ref{prop : secstr}, 
$$\langle R(\lambda_{\alpha},J\lambda_{\alpha})J\lambda_{\alpha},\lambda_{\alpha}\rangle\to -\infty$$
 as $n\to \infty$. By symmetries of the Riemann curvature tensor (see e.g. \cite[\S II.1]{chaveldg}) we have: 
\begin{eqnarray*}
\langle R(J\lambda_{\alpha},\lambda_{\alpha})\lambda_{\alpha},J\lambda_{\alpha}\rangle&=&\langle R(\lambda_{\alpha},J\lambda_{\alpha})J\lambda_{\alpha},\lambda_{\alpha}\rangle,\;\text{and}\\
 \langle R(J\lambda_{\alpha},\lambda_{\alpha})J\lambda_{\alpha},\lambda_{\alpha}\rangle&=&\langle R(\lambda_{\alpha},J\lambda_{\alpha})\lambda_{\alpha},J\lambda_{\alpha}\rangle=-\langle R(\lambda_{\alpha},J\lambda_{\alpha})J\lambda_{\alpha},\lambda_{\alpha}\rangle.
 \end{eqnarray*} 
 Then the four terms in the first bullet add up to 
 $$(a_{\alpha,n}e_{\alpha,n}-b_{\alpha,n}d_{\alpha,n})^{2}R(\lambda_{\alpha},J\lambda_{\alpha})J\lambda_{\alpha},\lambda_{\alpha}\rangle,$$
  Now since 
  $$\langle R(\lambda_{\alpha},J\lambda_{\alpha})J\lambda_{\alpha},\lambda_{\alpha}\rangle\to -\infty$$
   as $n\to \infty$, the sum of the four terms above is non-positive for all $n$ sufficiently large. 
   \medskip
 
\noindent{\bf Proof of the second bullet point:}  By symmetries of the Riemann curvature tensor any evaluation of Riemann curvature tensor with three or four $\lambda_{\alpha}$ or $J\lambda_{\alpha}$ vanishes. For example $\langle R(\lambda_{\alpha},J\lambda_{\alpha})\lambda_{\alpha},\lambda_{\alpha}\rangle=0$. So the terms which are multiple of them also vanish. The only evaluation of the Riemann curvature tensor with two $\lambda_{\alpha}$ and $J\lambda_{\alpha}$ which are not in the first bullet are $\langle R(J\lambda_{\alpha},J\lambda_{\alpha})\lambda_{\alpha},\lambda_{\alpha}\rangle$ and $\langle R(\lambda_{\alpha},\lambda_{\alpha})J\lambda_{\alpha},J\lambda_{\alpha}\rangle$ which are $0$. So their multiples also vanish. The rest of terms converge to $0$ as $n\to \infty$ by the bounds on the coefficients we established above and the convergence statement in Proposition \ref{prop : secstr} (see also the paragraph after Proposition \ref{prop : secstr} about the convention for the evaluations of the curvature tensor). Thus these terms have arbitrarily small absolute value.

Now using the symmetries of the Riemann curvature tensor and the bounds in the above two bullets, each term of the expansion of $\langle R(\dot{\zeta}_{n},J_{n})J_{n},\dot{\zeta}_{n} \rangle$ with one component $\lambda_{\alpha}$ or $J\lambda_{\alpha}$ ($\alpha\in \hat{\tau}$) is either non-positive or has arbitrarily small absolute value when $n$ is sufficiently large. The rest of the terms in the expansion of $\langle R(\dot{\zeta}_{n},J_{n})J_{n},\dot{\zeta}_{n} \rangle$ add up to $\langle R((\dot{\zeta_{n}})_{P},(J_{n})_{P})(J_{n})_{P},(\dot{\zeta}_{n})_{P} \rangle$. 
\medskip

We proceed to show that there is a $k_{3}<0$, such that 
\begin{equation}\label{eq : rcurv bd}\langle R((\dot{\zeta_{n}})_{P},(J_{n})_{P})(J_{n})_{P},(\dot{\zeta}_{n})_{P} \rangle\leq k_{3}\end{equation}
 for all $n$ sufficiently large. 
\medskip

We have $||J_{n}||^{2}\leq 1$, so $||(J_{n})_{P}||^{2}\leq 2$ by Lemma \ref{lem : asympprod}(\ref{asympprod : anorm}). By this inequality and (\ref{eq : c1}) we have  
$$\frac{\omega_{0}}{4}\leq ||(J_{n})_{P}||^{2}\leq 2.$$
 So after possibly passing to a subsequence $(J_{n})_{P}$ converge to a nonzero vector $v$ with $\frac{\omega_{0}}{4} \leq ||v||^{2}\leq 2$. 

Moreover $1\leq ||(\dot{\zeta}_{n})_{P}||^{2}\leq 2$ for $n$ sufficiently large by Corollary \ref{cor : geodP}. So the vectors $(\dot{\zeta}_{n})_{P}$ converge to a nonzero vector $w$ with $1\leq ||w||^{2}\leq 2$. 

Now we consider
\begin{eqnarray}\label{eq : Jz}
\langle (J_{n})_{N},(\dot{\zeta}_{n})_{N} \rangle&+&\langle (J_{n})_{P},(\dot{\zeta}_{n})_{N} \rangle\nonumber\\
&+&\langle (J_{n})_{N},(\dot{\zeta}_{n})_{P} \rangle+\langle (J_{n})_{P},(\dot{\zeta}_{n})_{P}\rangle= \langle J_{n},\dot{\zeta}_{n} \rangle\equiv 0.
\end{eqnarray} 
Since $\zeta_{n}(E')\subset U_{u_{n}}(\hat{\tau})$ and $u_{n}\to 0$, we have $a_{\alpha,n}\to 0$ and $b_{\alpha,n}\to 0$ by Lemma \ref{lem : geodnstr}. Therefore $||(\dot{\zeta}_{n})_{N}||^{2}\to 0$ as $n\to\infty$. Furthermore, since $||J_{n}||^{2}\leq 1$, for all $n$ sufficiently large $||(J_{n})_{N}||^{2}\leq \frac{1}{2}$ and $||(J_{n})_{P}||^{2}\leq \frac{1}{2}$ by Lemma \ref{lem : asympprod}(\ref{asympprod : anorm}). So by the Cauchy-Schwartz inequality, the first term and the second term of the sum in (\ref{eq : Jz}) go to $0$ as $n\to \infty$. We also have that $||J_{n}||^{2}\leq 1$, so $||(J_{n})_{P}||^{2}\leq 2$ and $||(J_{n})_{N}||^{2}\leq 2$ by Lemma \ref{lem : asympprod}(\ref{asympprod : anorm}). Then by Lemma \ref{lem : asympprod}(\ref{asympprod : aprod}) the third term of the sum in (\ref{eq : Jz}) goes to $0$ as $n\to \infty$. Having these bounds, from (\ref{eq : Jz}) we may conclude that 
$$\langle (J_{n})_{P},(\dot{\zeta}_{n})_{P} \rangle \to 0$$
 as $n\to \infty$. Therefore $\langle v,w \rangle=0$.
\medskip

As we saw above the evaluations of the Riemann curvature tensor 
$$\langle R((\dot{\zeta_{n}})_{P},(J_{n})_{P})(J_{n})_{P},(\dot{\zeta}_{n})_{P}\rangle$$
converge to $\langle R(v,w)w,v\rangle$ as $n\to\infty$. The vectors $v$ and $w$ are based at the point $p$ in the $u'-$thick part of the $\hat{\tau}-$stratum where all of the sectional curvatures of the stratum are bounded above by some $k_{2}<0$. Moreover, as we saw above, $||v||^{2}\leq 2$, $||w||^{2}\leq 2$ and $\langle v,w\rangle=0$, so 
$$|v\wedge w|^{2}=||v||^{2}||w||^{2}-\langle v,w \rangle^{2}\leq 4.$$
 Then 
$$\langle R(v,w)w,v\rangle=\kappa(v,w)|v\wedge w|^{2} \leq 4k_{2}.$$
Therefore, (\ref{eq : rcurv bd}) holds for $k_{3}=3k_{2}<0$. 
 
The two bullets above and (\ref{eq : rcurv bd}) show that for all $n$ sufficiently large, 
$$\langle R(\dot{\zeta}_{n},J_{n})J_{n},\dot{\zeta}_{n} \rangle$$ is the sum of terms which are either 
 negative or have arbitrarily small absolute value and the term 
 $$\langle R((\dot{\zeta_{n}})_{P},(J_{n})_{P})(J_{n})_{P},(\dot{\zeta}_{n})_{P} \rangle \leq k_{3}<0.$$ Thus for all $n$ sufficiently large, $\langle R(\dot{\zeta}_{n},J_{n})J_{n},\dot{\zeta}_{n} \rangle\leq \frac{k_{3}}{2}$. Moreover, $||\dot{\zeta}_{n}||^{2}\equiv 1$, $||J_{n}||^{2}\geq \frac{\omega_{0}}{2}$ and $\langle \dot{\zeta}_{n},J_{n}\rangle\equiv 0$, so $|\dot{\zeta}_{n}\wedge J_{n}|^{2}\geq \frac{\omega_{0}}{2}$. Thus 
 $$\kappa(\dot{\zeta}_{n}(s_{n}),J_{n}(s_{n}))\leq \frac{k_{3}}{\omega_{0}}<0.$$
 But this contradicts the assumption (\ref{eq : curv zeta J at sn}). The proof of the lemma is complete.
  \end{proof}

\begin{remark}
An explicit upper bound in terms of the systole of the Riemann surface representing the point in the Teichm\"{u}ller space where the curvature tensor and sectional curvatures are evaluated is not available; see the introduction of \cite{wpcurv}. So in the preceeding lemma we used compactness arguments to deduce the bounds on sectional curvatures.
\end{remark}

\noindent\textbf{Case 2:} We have $||J_{P}(s)||^{2}\leq \frac{\omega_{0}}{4}$ for some $s\in E$. 
\medskip

Let 
$$s_{0}=\min\{s : s\in E \;\text{and}\; ||J_{P}(s)||^{2}\leq \frac{\omega_{0}}{4}\}.$$
 Note that since $||J_{P}(0)||^{2}\geq \omega_{0}$, and $||J_{P}(s)||^{2}$ is a continuous function of $s$, we have $s_{0}>0$. We prove that the inequality (\ref{eq : Js<Jb}) holds using the following lemma:
 
 \begin{lem} \label{lem : intJ'}
There is a $\Phi>0$ depending only on $\omega_{0}$ and $b$ such that 
$$\int_{0}^{s}||J'||^{2}ds>\Phi$$ 
for every $s\geq s_{0}$.
\end{lem}
 The negative curvature of the WP metric and (\ref{eq : J''}) imply that 
$$\frac{1}{2}\frac{d^{2}}{ds^{2}}||J(s)||^{2}\geq ||J'(s)||^{2}.$$
 Integrating both sides of the above inequality on the interval $[0,s]$, we get 
 $$\frac{d}{ds}||J(s)||^{2}-\frac{d}{ds}||J||^{2}|_{s=0}\geq 2\int_{0}^{s}||J'||^{2} ds.$$
  By (\ref{eq : J'}), $\frac{d}{ds}||J||^{2}|_{s=0}=0$, so from the above inequality we obtain
  $$\frac{d}{ds}||J(s)||^{2}\geq 2\int_{0}^{s}||J'||^{2} ds.$$
 Moreover, for every $s\geq s_{0}$, $\int_{0}^{s}||J'||^{2} ds\geq\Phi$  by Lemma \ref{lem : intJ'}. Therefore 
 \begin{equation}\label{eq : der sqnormJ < 2phi}\frac{d}{ds}||J(s)||^{2}\geq 2\Phi.\end{equation}
 Now suppose that $b'\leq s_{0}$. Then integrating both sides of the iequality (\ref{eq : der sqnormJ < 2phi}) on the interval $[b',b]$ we get
 $$||J(b)||^{2}\geq||J(b')||^{2}+2\Phi (b-b').$$ 
Moreover, we have that $||J(b)||^{2}=1$, so 
$$1\geq ||J(b')||^{2}+2\Phi (b-b').$$
 Consequently, for $\delta=1-2\Phi (b-b')$, (\ref{eq : Js<Jb}) holds. Note that since $\omega_{0}<1$, we have $\Phi<\frac{1}{65b}$ (see the proof of Lemma \ref{lem : intJ'} for the value of $\Phi$), and then $\delta>0$. Furthermore, note that $\delta<1$ and that $\delta$ decreases as $b'\to 0$.
 \medskip

We finish Case 2 by proving Lemma \ref{lem : intJ'}.

\begin{proof}[Proof of Lemma \ref{lem : intJ'}]
By (\ref{eq : normJbd}) we have that $||J(s_{0})||^{2}\geq \omega_{0}$. Moreover by the setup of $s_{0}$ we have $||J_{P}(s_{0})||^{2}\leq \frac{\omega_{0}}{4}$. Furthermore, by (\ref{eq : normJbd}) we have $||J(s_{0})||^{2}\leq 1$. Then using Lemma \ref{lem : asympprod}(\ref{asympprod : anorm}) we may deduce that $||J_{N}(s_{0})||^{2}\leq 2$ and $||J_{P}(s_{0})||^{2}\leq 2$. Thus by Lemma \ref{lem : asympprod}(\ref{asympprod : aprod}), 
$$\langle J_{N}(s_{0}),J_{P}(s_{0})\rangle\to 0$$
 as $u\to 0$. From this limit we conclude that given $q>0$, for $\epsilon$ sufficiently small we have that
\begin{equation}\label{eq : JN.JP at s0} |\langle J_{N}(s_{0}),J_{P}(s_{0}) \rangle |\leq q.\end{equation}
Now consider the equality 
$$||J(s_{0})||^{2}=||J_{N}(s_{0})||^{2}+||J_{P}(s_{0})||^{2}+2\langle J_{N}(s_{0}),J_{P}(s_{0})\rangle.$$
 Then using the bounds we mentioned above we get
 \begin{equation}\label{eq : IJN(s0)I}||J_{N}(s_{0})||^{2}\geq \frac{3\omega_{0}}{4}-q.\end{equation}

It follows from the way we constructed the Jacobi field $J$ that $J(0)=||J(0)||\dot{g}|_{\pi(p)}$. Moreover, $g([0,T])\subset U_{\epsilon}(\sigma)$. Also by (\ref{eq : normJbd}), $||J(0)||^{2}\leq 1$. So by Lemma \ref{lem : geodnstr}  given $q>0$ for $\epsilon$ sufficiently small, 
\begin{equation}\label{eq : normJN at 0}||J_{N}(0)||^{2}\leq q.\end{equation} 

Now we have that
\begin{eqnarray}\label{eq : intJ'N}
\frac{\omega_{0}}{2}&\leq&\int_{0}^{s_{0}}\frac{d}{ds}||J_{N}||^{2}ds \leq \int_{0}^{s_{0}}2|\langle J_{N},J_{N}'\rangle|ds\\
&\leq&2\Big(\int_{0}^{s_{0}}||J_{N}'||^{2}ds\Big)^{1/2}\Big(\int_{0}^{s_{0}}||J_{N}||^{2}ds\Big)^{1/2}\nonumber\\
&\leq&2\Big(\int_{0}^{s_{0}}||J_{N}'||^{2}ds\Big)^{1/2}\sqrt{2b}.\nonumber
\end{eqnarray}
The first inequality above follows from the Fundamental theorem of calculus and the bounds $||J_{N}(s_{0})||^{2}\geq \frac{3\omega_{0}}{4}-q$ and $||J_{N}(0)||^{2}\leq q$ we established in (\ref{eq : JN.JP at s0}) and (\ref{eq : IJN(s0)I}), respectively, where $q=\frac{\omega_{0}}{8}$ and $u$ and $\epsilon$ are sufficiently small. The second inequality follows from the compatibility of the Riemannian metric and its Levi-Civita covariant derivative. The third inequality is the Cauchy-Schwartz inequality. For the last inequality note that $||J(s)||^{2}\leq 1$ for any $s\in[0,b]$ by (\ref{eq : normJbd}). Then Lemma \ref{lem : asympprod}(\ref{asympprod : anorm}) implies that for $u$ sufficiently small $||J_{N}(s)||^{2}\leq 2$. Furthermore, $s_{0}\leq b$. Thus 
$$\int_{0}^{s_{0}}||J_{N}||^{2}ds\leq 2b,$$
 giving us the last inequality.
\medskip

From the bound (\ref{eq : intJ'N}) we get
\begin{equation}\label{eq : intJN}\int_{0}^{s_{0}}||J_{N}'||^{2}ds\geq \frac{\omega_{0}^{2}}{32b}.\end{equation}
Moreover, by Lemma \ref{lem : asympprod}(\ref{asympprod : anorm}), we have that $||J'||^{2}\geq \frac{1}{2}||(J')_{N}||^{2}$. So we have
\begin{eqnarray}\label{eq : int|J'|^2}
\int_{0}^{s_{0}}||J'||^{2}ds&\geq& \frac{1}{2}\int_{0}^{s_{0}}||(J')_{N}||^{2}ds\nonumber\\
 &\geq& \frac{1}{2}\Big(\int_{0}^{s_{0}}||J'_{N}||^{2}ds-\int_{0}^{s_{0}} ||J'_{N}||^{2}-||(J')_{N}||^{2} ds\Big).
\end{eqnarray} 
By (\ref{eq : normJbd}), $||J||^{2}\leq 1$ along $\zeta$. Then by Lemma \ref{lem : N'} we have that 
$$\big| ||J'_{N}||^{2}-||(J')_{N}||^{2}\big|\to 0,$$
 as $u\to 0$. So the last integral in (\ref{eq : int|J'|^2}) goes to $0$ as $u\to 0$. Thus using the bound (\ref{eq : intJN}) from (\ref{eq : int|J'|^2}) we may deduce that for $u$ sufficiently small,
 \begin{equation}\label{eq : int|J'| 0 to s0}\int_{0}^{s_{0}}||J'||^{2}ds\geq\frac{\omega_{0}^{2}}{65b}.\end{equation}
Let $\Phi:=\frac{\omega_{0}^{2}}{65b}$. We have that $||J'||^{2}\geq 0$, hence for every $s\geq s_{0}$ we have 
$$\int_{0}^{s}||J'||^{2}ds\geq \int_{0}^{s_{0}}||J'||^{2}ds \geq\Phi.$$
where the second inequality follows from (\ref{eq : int|J'| 0 to s0}). This finishes the proof of the lemma.
 \end{proof}
 
Now we can conclude the proof of the theorem. Let $\delta$ be the maximum of the $\delta$'s we worked out in Cases 1 and 2, and $\frac{1+\omega_{0}}{2}$. Then the theorem follows from the conclusions of Cases 1 and 2 and the discussion before Case 1 about the case where $||J(0)||^{2}<\omega_0$. Note that $\delta$ depends only on $b',b$ and the constants $\omega_0$ and $k_0$. Moreover, the choice of $k_{0}$ depends only on $\bar{\epsilon}$ (coming from compactness arguments). Hence $\delta$ depends only on $b',b$ and $\bar{\epsilon}$.  
 \end{proof}

\section{Strongly asymptotic rays}\label{subsec : strongasymp}

In this section we prove Theorem \ref{thm : strongasym}. Throughout this section the pair $(\nu^{-},\nu^{+})$ is a narrow pair, where $\nu^{+}$ is a lamination as in $\S$\ref{subsec : endinv} and $\nu^{-}$ is a marking.  Our strategy to prove the theorem is the same as the one in \cite[\S 4]{bmm1} for the proof of the Recurrent Ending Lamination Theorem (Theorem \ref{thm : recel}). Let $r_{\nu^{\pm}}$ be a WP geodesic ray with prescribed itinerary as in Theorem \ref{thm : inftyray}, where the end invariant $(\nu^{-},\nu^{+})$ is narrow and has bounded annular coefficients. For simplicity denote $r_{\nu^{\pm}}$ by $r$. Suppose that $r'$ is a geodesic ray which is not strongly asymptotic to $r$ and its ending lamination is contained in $\nu^{+}$. In the Narrow Visibility Theorem (Theorem \ref{thm : narrowvis}) we show that there is a bi-infinite geodesic $g$ strongly asymptotic to the ray $r$ in the forward time and asymptotic to $r'$ in the backward time, i.e. $r$ is visible. Note that any measured lamination with bounded length along both $r$ and $r'$ has bounded length along $g$.  Let $\mathcal{L}$ be an ending measured lamination of $r'$. By Lemma \ref{lem : dlendmeas} $\mathcal{L}$ has bounded length along $r'$.  Now the forward ending lamination of $r'$ is $\nu^{+}$, so the support of $\mathcal{L}$ is a sublamination of $\nu^{+}$. Then Proposition \ref{prop : suppnu+bddl} implies that $\mathcal{L}$ has bounded length along $r$ too. Consequently $\mathcal{L}$ has bounded length along $g$ in both the forward and the backward time. But this violates the convexity of the length of measured laminations along WP geodesics (Theorem \ref{thm : convlf}). Thus we may conclude that the rays $r$ and $r'$ are strongly asymptotic. In other words, the forward ending lamination determines the strong asymptotic class of rays with narrow end invariant and bounded annular coefficients.

\begin{definition}
Let $\mathcal{X}$ be a metric space. Two geodesic rays $r,r':[0,\infty) \to \mathcal{X}$ parametrized by arc-length are asymptotic if there are positive $d$ and $T$ such that 
$$d(r(t),r'(t)) \leq d,$$
for every $t \geq T$. The rays $r$ and $r'$ are strongly asymptotic if for some $a\geq 0$, 
$$d(r(t+a),r'(t)) \to 0$$
 as $t \to \infty$.
\end{definition}

\begin{thm}\textnormal{(Asymptotic $\Longrightarrow$ strongly asymptotic)}\label{thm : asymsasym}

Let $r_{\nu^{\pm}}$ be a WP geodesic ray with prescribed itinerary with $A-$narrow end invariant and bounded annular coefficients. Then any WP geodesic ray $r'$ asymptotic to $r_{\nu^{\pm}}$ is also strongly asymptotic to $r_{\nu^{\pm}}$. 

\end{thm}

\begin{proof}

For simplicity of notation we denote $r_{\nu^{\pm}}$ by $r$. Assume that $r$ and $r'$ are parametrized by arc-length.

By Lemma \ref{lem : rvisitingU} given $A,d>0$ and  $\epsilon\leq \bar{\epsilon}$ (where $\bar{\epsilon}$ is the constant in Lemma \ref{lem : rvisitingU}) for any $k\in\mathbb{N}$ there is a time interval $[t_{k}^{-},t_{k}^{+}]$ of length at least $d$ and a large subsurface $Z_{k}$ such that $r([t^{-}_{k},t^{+}_{k}])\subseteq U_{\epsilon,\bar{\epsilon}}(\partial{Z_{k}})$. Moreover, $t_{k}^{+} <t_{k+1}^{-}$.

Let $n\in\mathbb{N}$. Consider the geodesic segments $[r(0),r'(n)]$ and $r([0,n])$ and let $U_{1,n}:\triangle_{1,n}\to \Teich(S)$ be the geodesic variation described in $\S$\ref{sec : ruledrec}. Similarly for the geodesic segments $[r(n),r'(n)]$ and $r([0,n])$ let $U_{2,n}:\triangle_{2,n}\to \Teich(S)$ be the geodesic variation described at the beginning of $\S$\ref{sec : ruledrec}. Let $\triangle_{n}=\triangle_{1,n}\cup\triangle_{2,n}$ and let $U_{n}:\triangle_{n}\to \Teich(S)$ be the map which restricts to $U_{1,n}$ on $\triangle_{1,n}$ and restricts to $U_{2,n}$ on $\triangle_{2,n}$; see Figure \ref{fig : ruledsurface}.

\begin{figure}
\centering
\scalebox{0.2}{\includegraphics{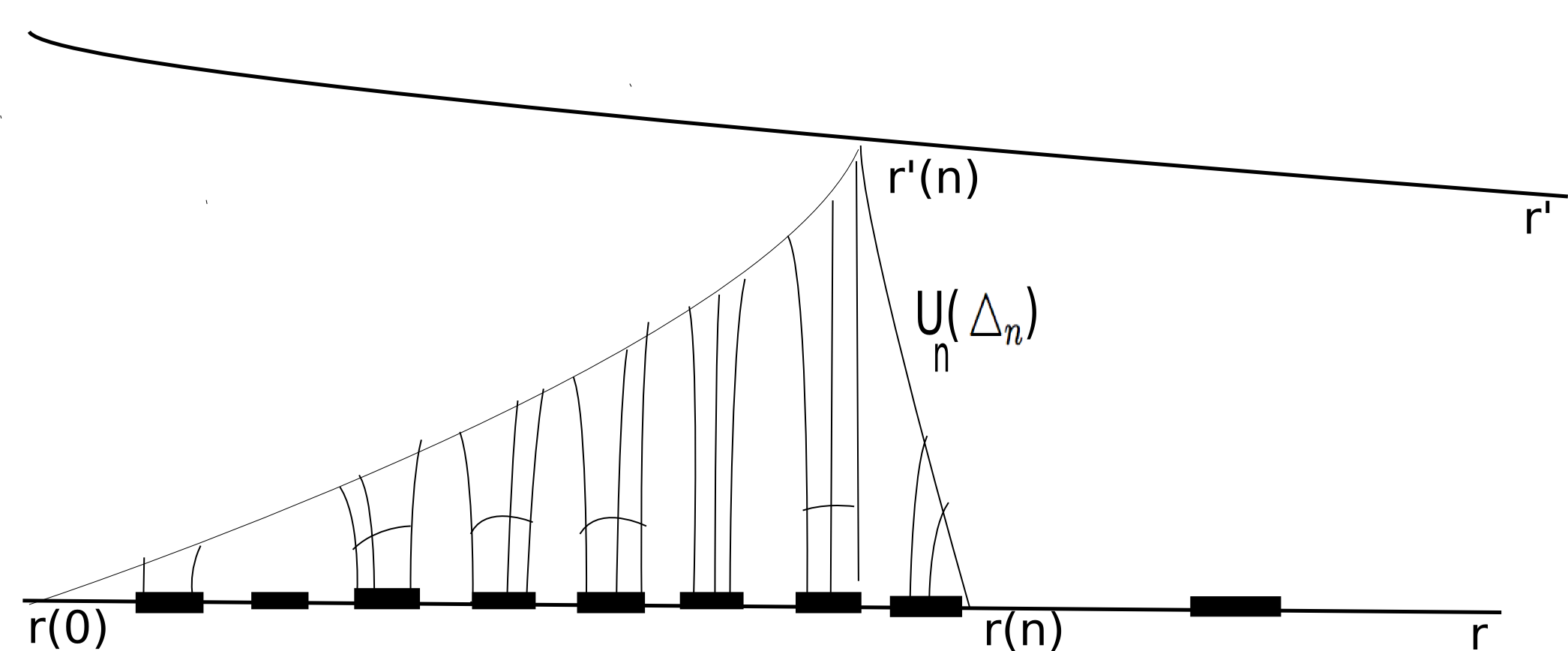}}
\caption{ The map $U_{n}:\triangle_{n}\to\Teich(S)$ is a variation of geodesics. The bold intervals are $r([t_{k}^{-},t_{k}^{+}])\subset U_{\epsilon,\bar{\epsilon}}(\partial{Z}_{k})$. The function $f_{n}:[0,T_{n}]\to\Teich(S)$ is a parametrization of  $[r(0),r'(n)]$. For $k,n$, the interval $V_{k,n}$ is so that for any $t\in V_{k,n}$, the nearest point to $f_{n}(t)$ is on $r([t_{k}^{-},t_{k}^{+}])$. Moreover, $W_{k,n}\equiv [0,b]$. The rectangular regions are $U_{n}(V_{k,n}\times W_{k,n})$.}
\label{fig : ruledsurface}
\end{figure}

Let $f_{n}:[0,T_{n}]\to \Teich(S)$ be a parametrization of the geodesic segment $[r(0),r'(n)]$ by arclength. Let $V_{k,n}\subseteq [0,T_{n}]$ be a subinterval such that $\pi(f_{n}(t))\in r([t_{k}^{-},t_{k}^{+}])$ for any $t\in V_{k,n}$.  Let $\lambda_{n}(t)$ be length of the interval $[f_{n}(t),\pi(f_{n}(t))]$. Let $k,n\in\mathbb{N}$ and $b>0$ be so that $b\leq \lambda_{n}(t)$ for every $t\in V_{k,n}$. Let $W_{k,n}\equiv [0,b]$. Moreover let the map 
$$u_{k,n}:=U_{n}|_{V_{k,n}\times W_{k,n}}:V_{k,n}\times W_{k,n}\to\Teich(S).$$
 We denote $V_{k,n}\times W_{k,n}$ with the pull back metric through $u_{k,n}$ by $\square_{k,n}$, and call $b$ the height of $\square_{k,n}$. Recall the measure $m$ introduced after Theorem \ref{thm : pullbackm} in $\S$\ref{sec : ruledrec}. Let $l_{b,k,n}$ be the length of $u_{k,n}(V_{k,n}\times b)$. By Theorem \ref{thm : pullbackm}(\ref{pmetric : distinc}) the length of $l_{k,n,s}$ is increasing in $s$ and is greater that $\length(V_{k,n}\times 0)$. Then since each interval $[t_{k}^{-},t_{k}^{+}]$ has length at least $d$, we may choose the interval $V_{k,n}$ so that 
 $$l_{k,n,b}=u_{k,n}(V_{k,n}\times b)\equiv d$$
  for any $k,n$ as above.

\begin{lem}\label{lem : K0}
 There is a $K_{0}>0$ depending on $b$ and $d$, so that for any $k,n$ as above 
 $$\int\int_{\square_{k,n}}-\kappa dA+\int_{V_{k,n}} m\geq K_{0}.$$
 \end{lem}
\begin{remark}
In \cite[\S 4]{bmm1} the recurrence to the thick part of the moduli space where all of the sectional curvatures are bounded above by a negative constant is used to produce regions with the above property. 
\end{remark}
\begin{proof}
For every $b'\in(0,b]$ consider the Jacobi field 
$$J_{b',k,n}=\frac{\partial u_{k,n}}{\partial t}\Big|_{u_{k,n}(t,b')},$$
 which is tangent to the path $u_{k,n}(.,b')$. For each $t\in V_{k,n}$ consider the geodesic $u_{k,n}(t,.)$. Let $\delta$ be the constants from Theorem \ref{thm : normJbd} for the geodesic segments $u_{k,n}(t,.)$. Note that since $V_{k,n}$ is a compact interval, we may choose $b>0$ uniform over $V_{k,n}$ such that Theorem \ref{thm : normJbd} holds for $b$ and any $b'\leq b$. Also given $b'$ we may choose $\delta$ uniform over $V_{k,n}$. Then by the theorem, 
$$||J_{b',k,n}||^{2}\leq \delta||J_{b,k,n}||^{2}.$$
Moreover by Theorem  \ref{thm : normJbd}, $\delta$ decreases as $b'\to 0$, so we may choose $\delta$ uniform for all $b'$. Integrating both sides of the above inequality over the interval $V_{k,n}$ and taking square root we get 
$$l_{b',k,n}\leq \sqrt{\delta}  l_{b,k,n}.$$ 
Now as we defined before Lemma \ref{lem : lengthcurv}, let $l_{0,k,n}=\lim_{s\to 0}l_{s,k,n}$. Then taking limit as $b'\to 0$ of both sides of the above inequality we obtain 
$$l_{0,k,n}\leq \sqrt{\delta} l_{b,k,n}.$$
 Subtracting both sides of this inequality from $l_{b,k,n}$ we get 
$$l_{b,k,n}-l_{0,k,n}\geq l_{b,k,n}-\sqrt{\delta}l_{b,k,n}.$$
Then since  $l_{b,k,n}\equiv d$ for all $k,n$, from the above inequality we obtain 
$$l_{b,k,n}-l_{0,k,n}\geq (1-\sqrt{\delta})d.$$
Then Lemma \ref{lem : lengthcurv} guarantees that the lemma holds for $K_{0}=\frac{(1-\sqrt{\delta})d}{b}$.
\end{proof}

We proceed to finish the proof of the theorem by contradiction. Suppose that the geodesic rays $r$ and $r'$ are not strongly asymptotic. Then since the distance function in a $\CAT(0)$ space is a convex function, there are $b_{0}>0$ and $T>0$ such that 
$$d_{\WP}(r(t),r'(t))>b_{0}$$ 
for all $t\geq T$. Let $b=\frac{b_{0}}{2}$. For any $n\in\mathbb{N}$, let $N_{n}$ be the number of ruled rectangles $\square_{k,n} \subset \triangle_{n}$ with hight $b$. Then by Lemma \ref{lem : K0} we have
$$N_{n} K_{0}\leq \sum_{k=1}^{N_{n}}(\int\int_{\square_{k,n}}-\kappa dA+\int_{V_{k,n}}m).$$
 Now note that $-\kappa$ is a positive function and $m$ is a positive measure, and moreover $\bigcup_{k}\square_{k,n}\subset \triangle_{n}$. Thus, we have
$$\sum_{k=1}^{N_{n}}(\int\int_{\square_{k,n}}-\kappa dA+\int_{V_{k,n}}m)\leq \int \int_{\triangle_{n}}-\kappa dA+\int_{[0,n]} m.$$
By the formula (\ref{eq : gb}) the right-hand side is bounded above by $\pi$ independent of $n$. Thus $N_{n}\leq \frac{\pi}{K_{0}}$ for all $n$. On the other hand, the $\CAT(0)$ comparison for $\triangle_{n}$ shows that $N_{n}$ would be made arbitrarily large by increasing $n$. But this contradicts the upper bound for $N_{n}$ we obtained above. The fact that $r$ and $r'$ are strongly asymptotic follows from this contradiction. The proof of the theorem is complete.
 \end{proof}

\begin{thm}\textnormal{(Narrow visibility)}\label{thm : narrowvis}

Let $r_{\nu^{\pm}}$ be a WP geodesic ray with prescribed itinerary, where the end invariant $(\nu^{-},\nu^{+})$ is narrow and has bounded annular coefficients. Let $r'$ be a WP geodesic ray which is not strongly asymptotic to $r$. Then there is a bi-infinite geodesic in $\overline{\Teich(S)}$ which is strongly asymptotic to $r_{\nu^{\pm}}$ in the forward time and asymptotic to $r'$ in the backward time. In other words, $r_{\nu^{\pm}}$ is visible.

We know that the restriction of $\nu^{+}$ to a large subsurface $Y$ is in $\mathcal{EL}(Y)$.  Suppose that a measured lamination $\mathcal{L}\in\mathcal{ML}(Y)$ has bounded length along both $r_{\nu^{\pm}}$ and $r'$. Then $\mathcal{L}$ has bounded length along $g$.
\end{thm}
 
\begin{proof}

For simplicity we denote $r_{\nu^{\pm}}$ by $r$. Let $r$ and $r'$ both be parametrized by arc-length. We may omit finitely many of intervals $[t^{-}_{k},t^{+}_{k}]$ and assume that there is $b_{0}>0$ such that for each $k$ the distance between the geodesic segment $r([t^{-}_{k},t^{+}_{k}])$ and $r'$ is at least $b_{0}$. The reason is that in a $\CAT(0)$ space the distance between any two geodesic rays is a convex function, so the distance between $r$ and $r'$ is a monotonic function after some time. Then since the rays $r$ and $r'$ are not strongly asymptotic, the distance function is bounded below.

We briefly recall the setup from Theorem \ref{thm : asymsasym}. For any $n\in\mathbb{N}$, consider the geodesic segments $[r(0),r'(n)]$ and $r([0,n])$, and the variation of geodesics obtained from geodesic segments connecting a point on $[r(0),r'(n)]$ to its nearest point on $r([0,n])$. Similarly, consider the geodesic segments $r([0,n])$ and $[r(n),r'(n)]$ and the variation of geodesics obtained from the nearest point projection. Let $U_{n}:\triangle_{n}\to \mathcal{X}$ be the map defined by putting together the two geodesic variations. Let $f_{n}:[0,T_{n}]\to \Teich(S)$ be a parametrization of $[r(0),r'(n)]$ by arclength. For each $k\in\mathbb{N}$ let $V_{k,n}\subseteq [0,T_{n}]$ be such that $\pi(f_{n}(t))\in r([t_{k}^{-},t_{k}^{+}])$ for any $t\in V_{k,n}$. Suppose that there is a $b>0$ so that $b\leq \lambda_{n}(t)$ for all $t\in V_{k,n}$. Then let $W_{k,n} \equiv [0,b]$. We may assume that the length of all intervals $V_{k,n}$ is $d$. Let 
$$u_{k,n}:=U_{n}|_{V_{k,n}\times W_{k,n}}:V_{k,n}\times W_{k,n}\to \Teich(S).$$
 We denote $V_{k,n}\times W_{k,n}$ with the metric pull back through $u_{k,n}$ by $\square_{k,n}$. Let $g_{n}:[a_{n},b_{n}]\to\Teich(S)$ be a parametrization by arc-length of $[r(n),r'(n)]$.

First we show that after possibly passing to a subsequence there are parameters $\hat{t}_{n}\in[a_{n},b_{n}]$ so that the points $g_{n}(\hat{t}_{n})$ converge to a point $z \in\overline{\Teich(S)}$.

Let $T>0$ be such that $d(r(t),r'(t))>b_{0}$ for every $t>T$. Let $b=\frac{b_{0}}{2}$. We claim that there is a $\hat{k}\in\mathbb{N}$ such that $u_{\hat{k},n}(\square_{\hat{k},n}) \cap g_{n} \neq \emptyset$ for all sufficiently large $n$.  

Let $K_{0}>0$ be the constant in Lemma \ref{lem : K0} which we proved in the course of the proof of Theorem \ref{thm : asymsasym}. Let $\hat{k},n\in\mathbb{N}$ be so that $g_{n} \cap u_{\hat{k},n}(\square_{\hat{k},n})=\emptyset$. Let $N_{n}$ be the number of ruled rectangles $\square_{k,n}\subset \triangle_{n}$ with hight $b$. Then since $\bigcup_{k}\square_{k,n}\subset\square_{n}$, as we saw in the proof of Theorem \ref{thm : asymsasym}, we have
$$N_{n} K_{0}\leq \int\int_{\triangle_{n}} -\kappa dA+\int_{[T,n]}m.$$
By the formula (\ref{eq : gb}) the right-hand side above is bounded above by $\pi $ independent of $n$. Thus $N_{n}\leq \frac{\pi}{K_{0}}$ for all $n$ sufficiently large. On the other hand, by the $\CAT(0)$ comparison for $\triangle_{n}$, $N_{n}$ can be made arbitrarily large by increasing $n$.  We may also assume that $n>\hat{k}$. So for $\hat{k}$ sufficiently large we get a contradiction to the upper bound for $N_n$ above, and the claim follows.
\medskip

As we saw above there is an integer $\hat{k}$ so that for all $n\in\mathbb{N}$ sufficiently large, the geodesic segment $g_{n}$ and $\square_{\hat{k},n}$ intersect.  Let the point $z_{n}\in\square_{\hat{k},n}$ and the parameter $\hat{t}_{n}$ be so that $z_{n}=g_{n}(\hat{t}_{n})$. Then since the hight of each $\square_{\hat{k},n}$ is less than $b$, all of the points $z_{n}$ are in the $b-$neighborhood of $r([t^{-}_{\hat{k}},t^{+}_{\hat{k}}])$ which is a compact subset of $\overline{\Teich(S)}$. Therefore after possibly passing to a subsequence we may assume that the points $z_{n}$ converge to a point $z$ in the $b-$neighborhood of $r([t^{-}_{\hat{k}},t^{+}_{\hat{k}}])$. 
\medskip

Now we proceed to show that the geodesic segments $g_{n}$ converge to a geodesic $g$ passing through $z$. We reparametrize $g_{n}$ by arc-length such that $g_{n}(0)=z_{n}$. Let $t \in [0,\infty)$ and let the integer $N_{t}\geq 1$ be such that $d(z,r(n))>t$ for each integer $n\geq N_{t}$. We show that the sequence $\{g_{n}(t)\}_{n=N_{t}}^{\infty}$ is a Cauchy sequence. To see this, let $h^{+}_{n}$ be the geodesic segment joining $z$ to $r(n)$ parametrized by arc-length. As is shown in Lemma 8.3 in \cite[\S II.8]{bhnpc} $\{h^{+}_{n}(t)\}_{n=N_{t}}^{\infty}$ is a Cauchy sequence. Furthermore, for each $n\geq N_{t}$, let $g^{+}_{n}=g_{n}|_{[0,b_{n}]}$ be the parametrization of the geodesic segment $[z_{n},r(n)]$. Then by the $\CAT(0)$ comparison for the triangle with vertices $z$, $ z_{n}$ and $r(n)$, 
$$d(h^{+}_{n}(t),g^{+}_{n}(t))< d(z_{n},z).$$
 Then since $z_{n}\to z$, we have that $\{g_{n}(t)\}_{n=N_{t}}^{\infty}$ is a Cauchy sequence.
Let $t\in(-\infty,0]$ and let the integer $N_{t}\geq 1$ be such that $d(z,r(n))>-t$ for each integer $n\geq N_{t}$. A similar argument to the one above shows that $\{g_{n}(t)\}_{n=N_{t}}^{\infty}$ is a Cauchy sequence. Therefore for each $t\in\mathbb{R}$ the sequence $g_{n}(t)$ is convergent in $\overline{\Teich(S)}$. The point-wise limit of a sequence of parametrized geodesics in a complete $\CAT(0)$ space is a parametrized geodesic. Moreover the convergence is as parametrized geodesics. Thus the geodesics $g_{n}$ converge to a bi-infinite geodesic $g:\mathbb{R}\to \overline{\Teich(S)}$ as parametrized geodesics. See the proof of Proposition 8.2 in \cite[\S II.8]{bhnpc} for more detail. The geodesic $g$ is the limit of the geodesic segments $[r(n),r'(n)]$, so $g^{+}=g|_{[0,\infty)}$ is asymptotic to $r$ and $g^{-}=g|_{[0,-\infty)}$ is asymptotic to $r'$. Moreover, by Theorem \ref{thm : asymsasym}, $g^{+}$ is strongly asymptotic to $r$. The first assertion of the theorem is proved.
\medskip

We continue by proving the second assertion of the theorem. First, suppose that the subsurface $Y$ is $S$. In the discussion above for the proof of the first assertion of the theorem, we may choose $b$ small enough that the closure of the $b-$neighborhood of the $\bar{\epsilon}-$thick part of the Teichm\"{u}ller space does not intersect any completion stratum. Each point $z_{n}$ is in the thick part of the Teichm\"{u}ller space, and therefore $z\in\Teich(S)$. Thus the measured lamination $\mathcal{L}$ has bounded length at $z$. By the assumption of the theorem $\mathcal{L}$ has uniformly bounded length at the points $r(n)$. So by Theorem \ref{thm : convlf} (Convexity of length-functions) $\mathcal{L}$ has uniformly bounded length along the geodesic segments $[z,r(n)]$. Then by continuity of length-functions $\mathcal{L}$ has bounded length along $g^{+}$. A similar argument shows that $\mathcal{L}$ has bounded length along $g^{-}$. Therefore $\mathcal{L}$ has bounded length along $g$.

Now suppose that $Y$ is a proper subsurface of $S$. Let $\rho$ be a hierarchy path between the narrow pair $(\nu^{-},\nu^{+})$. There is an $N\in\mathbb{N}$ such that $\partial{Y}\subset\rho(i+N)$ for all $i\in\mathbb{N}$; see \cite[\S 5]{mm2}. Then Theorem \ref{thm : inftyray} guarantees that for any $k$ sufficiently large, $r([t_{k}^{-},t_{k}^{+}])\subset U_{\epsilon,\bar{\epsilon}}(\partial{Z})$ for some subsurface $Z$ with $\partial{Y}\subseteq\partial{Z}$. In the discussion for the proof of the first assertion of the theorem, we may choose $b$ small enough that the $b-$neighborhood of $U_{\epsilon,\bar{\epsilon}}(\partial{Z})$ only intersects $\sigma-$strata where $\sigma\subseteq \partial{Z}$. Then choosing $\hat{k}$ large enough each $z_{n}$ is in the $b-$neighborhood of $r([t_{\hat{k}}^{-},t_{\hat{k}}^{+}])\subset U_{\epsilon,\bar{\epsilon}}(\partial{Z})$. Thus $z$ the limit of the points $z_{n}$ after possibly passing to a subsequence is a point in a $\sigma-$stratum where $\sigma\subseteq\partial{Z}$ (note that $\sigma$ could be empty). Then the non-refraction property of completion strata (Theorem 5.2 in \cite[Chapter 5]{wol}) guarantees that $g\subset \mathcal{S}(\sigma)$. 

We claim that $\sigma\subseteq\partial{Y}$. To see this, note that the geodesic segments $[r(n),r'(n)]$ after possibly passing to a subsequence converge to $g$. Thus the points $r(n)$ converge to a point in $\mathcal{S}(\sigma)$. Now suppose that there is a curve $\alpha\in\sigma$ with $\alpha\notin\partial{Y}$. First, $\alpha$ is a pinching curve of $r$. Thus by Definition \ref{def : endlam}, $\alpha$ is a component of the forward ending lamination of $r$. Second, since $\sigma$ and $\partial{Y}$ are subsets of $\partial{Z}$, $\alpha$ does not overlap $\partial{Y}$. Moreover, $Y$ is a large subsurface. Thus we have $\alpha\subset Y$. But this contradicts our assumption that the restriction of the forward ending lamination $\nu^{+}$ of $r$ to $Y$ is minimal filling. To see this, note that $\xi(Y)>1$, so $\alpha$ does not fill $Y$. Therefore $\alpha$ is a component of $\nu'$ the restriction of $\nu^{+}$ to $Y$ which contradicts the fact that $\nu'$ is minimal. The claim follows. 
\medskip

Now since the support of $\mathcal{L}$ is the lamination $\nu'$ in the subsurface $Y$, the above claim guarantees that none of the pinched curves at $z$ (the curves in $\sigma$) intersects $\mathcal{L}$. This implies that $\mathcal{L}$ has bounded length at $z$. Then an argument similar to what we gave above, where we assumed that $z$ is in the thick part of the Teichm\"{u}ller space, implies that $\mathcal{L}$ has bounded length along $g$. The second assertion of the theorem is proved.
\end{proof}

\begin{prop}\label{prop : suppnu+bddl}
Let $r_{\nu^{\pm}}$ be a ray with prescribed itinerary where the end invariant $(\nu^{-},\nu^{+})$ is narrow and has bounded annular coefficients. Any measured lamination supported on a sublamination of $\nu^{+}$ has bounded length along $r_{\nu^{\pm}}$. 
\end{prop}

\begin{proof} As in the proof of Lemma 4.5 in \cite{bmm1} we consider the ray $\bar{r}$ as follows: since $(\nu^{-},\nu^{+})$ is a narrow pair, $\nu^{+}$ has a minimal component $\nu'$ that fills a large subsurface $Y$. Let $\Sigma_{\nu'}$ be the simplex of the projective classes of measures supported on $\nu'$ in $\mathcal{PML}(Y)$. Let $\bar{\mathcal{L}} \in \mathcal{ML}(Y)$ be a representative of the projective class determined by a point in the interior of the top-dimensional face of $\Sigma_{\nu'}$. Then $\bar{\mathcal{L}}$ is a linear combination of all ergodic measures supported on $\nu'$ with  positive coefficients. Let $\{\gamma_{n}\}_{n}\subset\mathcal{C}_{0}(Y)$ be a sequence of simple closed curves so that the projective classes $[\gamma_{n}]$ converge to $[\bar{\mathcal{L}}]$. Let $c_{n}$ be a hyperbolic surface pinched at a pants decomposition that contains $\gamma_{n}\cup \partial{Y}$. Let $r_{n}$ be a parametrization of $[x,c_{n}]$ by arclength, and let $\bar{r}$ be a limit of the geodesic segments $r_{n}$ after possibly passing to a subsequence in the visual sphere of the Teichm\"{u}ller space at $x$. By Lemma \ref{lem : inftyray} $\bar{r}$ is an infinite ray and $\bar{\mathcal{L}}$ has bounded length along $\bar{r}$. This implies that any ergodic measure supported on $\nu'$ has bounded length along $\bar{r}$. Any measure supported on $\nu'$ is a linear combination of ergodic measures of $\nu'$ with non-negative coefficients. So we may conclude that any measure supported on $\nu'$ has bounded length along $\bar{r}$. Moreover by Lemma \ref{lem : inftyray} every $\alpha\in\partial{Y}$ has bounded length along $\bar{r}$. Any  sublamination of $\nu^{+}$ is the union of $\nu'$ and some of the curves in $\partial{Y}$. Therefore, any measure supported on a sublamination of $\nu^{+}$ has bounded length along $\bar{r}$.

 We proceed to show that $r_{\nu^{\pm}}$ and $\bar{r}$ are strongly asymptotic rays. Denote $r_{\nu^{\pm}}$ by $r$. If $r$ is not strongly asymptotic to $\bar{r}$, then by Theorem \ref{thm : narrowvis} there is a bi-infinite WP geodesic $g$ strongly asymptotic to $r$ in the forward time and asymptotic to $\bar{r}$ in the backward time. Let $\mathcal{L}$ be any ending measured lamination of $r$. By Lemma \ref{lem : dlendmeas}, $\mathcal{L}$ has bounded length along $r$. Moreover, $\mathcal{L}$ is supported on a sublamination of the forward ending lamination of $r$. By Lemma \ref{lem : inftyraynu} the forward lamination of $r$ contains $\nu'$, and therefore  is a sublamination of  $\nu^{+}$. Thus as we saw in the previous paragraph $\mathcal{L}$ has bounded length along $\bar{r}$. Then by the second assertion of Theorem \ref{thm : narrowvis}, $\mathcal{L}$ has bounded length along $g$. But this contradicts the convexity of length-functions along WP geodesics (Theorem \ref{thm : convlf}). This contradiction shows that $r$ and $\bar{r}$ are in fact strongly asymptotic. Now since both $r$ and $\bar{r}$ start at $x$, the convexity of the distance function between two geodesics in a $\CAT(0)$ space implies that $r=\bar{r}$.

 As we saw in the first paragraph of the proof any measured lamination supported on a sublamination of $\nu^{+}$ has bounded length along $\bar{r}$. By the second paragraph $r_{\nu^{\pm}}=\bar{r}$, thus any measured lamination supported on a sublamination of $\nu^{+}$ has bounded length along $r$, as was desired.
\end{proof}

Here we rephrase Theorem \ref{thm : strongasym} and prove it.
\begin{thmnarrowend}
Let $(\nu^{-},\nu^{+})$ be a narrow pair with bounded annular coefficients. Any WP geodesic ray $r'$ with forward ending lamination $\nu^{+}$ is strongly asymptotic to $r_{\nu^{\pm}}$.
\end{thmnarrowend}

\begin{proof} Denote $r_{\nu^{\pm}}$ by $r$. The lamination $\nu^{+}$ has a minimal component $\nu'$ that fills a large subsurface $Y$. We show that any other infinite ray $r'$ with forward ending lamination containing $\nu'$ is strongly asymptotic to $r$. Suppose not. Then by Theorem \ref{thm : narrowvis} there is a bi-infinite geodesic $g$ strongly asymptotic to $r$ in the forward time and asymptotic to $r'$ in the backward time. Let $\mathcal{L}$ be an ending measured lamination of $r'$. Then by Lemma \ref{lem : dlendmeas} the length of $\mathcal{L}$ is bounded along $r'$. Moreover $\mathcal{L}$ is supported on a sublamination of the forward ending lamination of $r'$ which is  $\nu^{+}$. Then Proposition \ref{prop : suppnu+bddl} implies that $\mathcal{L}$ has bounded length along $r$ as well. Then by the second assertion of Theorem \ref{thm : narrowvis}, $\mathcal{L}$ has bounded length along the bi-infinite geodesic $g$. But this contradicts convexity of the $\mathcal{L}-$length-function along $g$ (Theorem \ref{thm : convlf}). Thus we conclude that $r$ and $r'$ are strongly asymptotic.
\end{proof}

As an application of Theorem \ref{thm : strongasym} we provide a symbolic condition in terms of subsurface coefficients for divergence of WP geodesic rays in the moduli space. A geodesic ray is divergent if eventually leaves any compact subset of the moduli space. Recall that a pair $(\nu^{-},\nu^{+})$ is $A-$narrow if $d_{Z}(\nu^{-},\nu^{+})>A$ implies that $Z$ is a large subsurface.
\begin{thm}\textnormal{(Divergence condition)}
Given $A,R,R'>0$. Let $(\nu^{-},\nu^{+})$ be an $A-$narrow pair on a Riemann surface $S$ with $R'-$bounded annular coefficients and suppose that $d_{S}(\nu^{-},\nu^{+})\leq R$. Then a WP geodesic ray in $\mathcal{M}(S)$ with end invariant $(\nu^{-},\nu^{+})$ is divergent.
\end{thm}

\begin{proof}
Let $r_{\nu^{\pm}}$ be the infinite ray with prescribed itinerary and end invariant $(\nu^{-},\nu^{+})$. For simplicity we denote $r_{\nu^{\pm}}$ by $r$. By the narrow condition the restriction of $\nu^{+}$ to a large subsurface $Y$ is minimal and fills $Y$. Moreover the condition $d_{S}(\nu^{-},\nu^{+})\leq R$ guarantees that $Y$ is a proper subsurface of $S$. We proceed to show that given $\epsilon>0$ there is a time $T>0$ such that $\ell_{\alpha}(r(T))\leq\epsilon$ for every $\alpha\in \partial{Y}$. Let $\rho:[0,\infty]\to P(S)$ be a hierarchy path between $\nu^{-}$ and $\nu^{+}$. There is an $N\in\mathbb{N}$ such that for any component domain $Z$ of $\rho$ satisfying $J_{Z}\subseteq [N,\infty]$ we have $\partial{Y}\subseteq\partial{Z}$; see \cite[\S 5]{mm2}. 

Let $\epsilon>0$ and set $d=1$. It follows from Lemma \ref{lem : rvisitingU} that there are $L'>0$ and $\bar{w}>0$ with the property that for any subinterval $[m',n']\subset [0,\infty]$ with $n'-m'>L'$ there are $Z$ a component domain of $\rho$ and a time $T>0$ so that $\ell_{\alpha}(r(T))\leq\epsilon$ for all $\alpha\in\partial{Z}$. Furthermore, as we saw above if $m'>N$, then $\partial{Y}\subseteq \partial{Z}$. Thus for any $\alpha\in\partial{Y}$, $\ell_{\alpha}(r(T))\leq\epsilon$. For each $\alpha\in \partial{Y}$, the function $\ell_{\alpha}(r(t)):\mathbb{R}^{\geq 0} \to \mathbb{R}^{\geq 0}$ is a convex function, moreover by Lemma \ref{lem : inftyray}, $\ell_{\alpha}$ is bounded, thus it is a decreasing function. Therefore, 
$$\ell_{\alpha}(r(t))\leq \epsilon$$
 for all $t\in [T,\infty)$. Now choosing $\epsilon$ arbitrarily small, we may conclude that the projection $\hat{r}$ of the geodesic ray $r$ to $\mathcal{M}(S)$ is a divergent WP geodesic ray. 

Now by Theorem \ref{thm : strongasym} any geodesic ray $r'$ with the forward end invariant $\nu^{+}$ is strongly asymptotic to $r$. As we saw above $\hat{r}$ is divergent in $\mathcal{M}(S)$, which implies that the projection $\hat{r'}$ of $r'$ to $\mathcal{M}(S)$ is divergent as well. 
\end{proof}
On any surface $S$ with $\xi(S)>1$, in \cite[\S7.2]{wpbehavior} we constructed pairs $(\nu^{-},\nu^{+})$ where $\nu^{-}$ is a marking and  $\nu^{+}$ is a minimal filling lamination with the following properties:
\begin{itemize}
\item There is a list of large subsurfaces $\{Z_{i}\}_{i=1}^{\infty}$ so that for each $i$, $Z_{i}=Z_{i-1}\cap Z_{i+1}$, and $d_{Z_{i}}(\nu^{-},\nu^{+})\to \infty$ as $i\to \infty$. 
\item The remaining subsurface coefficients, including the annular ones, are uniformly bounded above. 
\end{itemize}
Furthermore, in $\S8$ of \cite{wpbehavior} we proved that there is a WP geodesic ray $r_{\nu^{\pm}}$ with end invariant $(\nu^{-},\nu^{+})$ as above and prescribed itinerary so that the projection $\hat{r}_{\nu^{\pm}}$  of $r_{\nu^{\pm}}$ to the moduli space is divergent. Here we show that any WP geodesic with end invariant $(\nu^{-},\nu^{+})$ is divergent in the moduli space.

\begin{thm}
Any WP geodesic ray with end invariant $(\nu^{-},\nu^{+})$ constructed in \cite[\S7.2]{wpbehavior} is divergent in $\mathcal{M}(S)$. 
\end{thm}
\begin{proof}
Since the only subsurfaces in which the pair $(\nu^{-},\nu^{+})$ has a big projection distance are large subsurfaces, the pair is narrow. Moreover the pair has bounded annular coefficients. Thus, by Theorem \ref{thm : strongasym} any other geodesic ray $r'$ with forward ending lamination $\nu^{+}$ is strongly asymptotic to $r_{\nu^{\pm}}$. Now since $\hat{r}_{\nu^{\pm}}$ is divergent in $\mathcal{M}(S)$, the geodesic ray $\hat{r}'$ is divergent as well. 
\end{proof}

\bibliographystyle{amsalpha}
\bibliography{reference}

\end{document}